\documentclass[a4paper]{article}
\usepackage[utf8]{inputenc}
\usepackage[a4paper, total={6.5in,9in}]{geometry}

\usepackage{natbib}

\usepackage{amssymb}
\usepackage{graphicx} 
\usepackage{xspace}
\usepackage{multirow}
\usepackage{bbm}
\usepackage{tabu}

\usepackage{algorithm}
\usepackage{algpseudocode}

\usepackage{graphicx}
\usepackage{subcaption}
\usepackage{enumitem}

\usepackage{multirow}
\usepackage{adjustbox}

\usepackage[colorinlistoftodos]{todonotes}

\usepackage{bm}
\usepackage{verbatim}
\usepackage{authblk}

\usepackage{amsmath}
\usepackage{hyperref}

\usepackage{lipsum}

\newcommand\blfootnote[1]{%
    \begin{NoHyper}
  \renewcommand\thefootnote{}\footnote{#1}%
  \addtocounter{footnote}{-1}%
    \end{NoHyper}
  }

\title{Weak form Shallow Ice Approximation models with an improved time step restriction}

\author[1,2]{Igor Tominec}
\author[1,2,3]{Josefin Ahlkrona}

\affil[1]{Stockholm University, Department of Mathematics, Stockholm, Sweden}
\affil[2]{Bolin Centre for Climate Research, Stockholm University, Stockholm, Sweden}
\affil[3]{Swedish e-Science Research Centre, Sweden}

\begin{document}
\maketitle
\abstract{
The Shallow Ice Approximation (SIA) model on strong form is commonly used for inferring
the flow dynamics of grounded ice sheets.
The solution to the SIA model is a closed-form expression for the velocity field.
When that velocity field is used to advance the ice surface in time,
the time steps have to take small values due to quadratic scaling in terms of the horizontal mesh size.
In this paper we write the SIA model on weak form, and add in the Free Surface Stabilization Algorithm (FSSA) terms.
We find numerically that the time step restriction scaling is improved from quadratic to linear,
but only for large horizontal mesh sizes. We then extend the weak form by adding the initially neglected normal stress terms.
This allows for a linear time step restriction across the whole range of the horizontal mesh sizes, leading
to an improved efficiency.
Theoretical analysis demonstrates that the inclusion of FSSA stabilization terms transitions the explicit time stepping treatment
of second derivative surface terms to an implicit approach.
Moreover, a computational cost analysis, combined with numerical results on stability and accuracy,
advocates for preferring the SIA models written on weak form over the standard SIA model.
}

\blfootnote{E-mail addresses:} 
\blfootnote{igor.tominec@math.su.se (Igor Tominec)} 
\blfootnote{ahlkrona@math.su.se (Josefin Ahlkrona)} 

\section{Introduction}

The Shallow Ice Approximation (SIA) problem is a commonly used momentum balance model which desribes the 
non-Newtonian, viscous, gravity driven flow of the ice in grounded ice sheets 
\citep{Hutter1983TheoreticalGM}. 
It is typically used either as a standalone model or in combined with the Shallow Shelf Approximation (SSA) in hybrid models for sea-level rise predictions \citep{tc-14-3071-2020,tc-14-3033-2020} on time-scales of a few hunded years as well as for paleoclimate spin-up simulations \citep{tc-13-1441-2019} and paleosimulations of 10 000 \citep{Weber2021} five million years \citep{Pollard2009}.
The SIA model is a simplification of the nonlinear (full) Stokes problem on the premise that an ice sheet is thin, 
neglecting all stress-components except vertical shear stresses.
Advantages of the SIA problem over the nonlinear full Stokes problem are 
that the SIA problem is linear and computationally less expensive to solve. Some of the disadvantages 
when compared to the 
nonlinear Stokes problem are: 
(i) the degraded model accuracy, (ii) that when coupled to the free-surface equation the simulation time steps have to be taken very small at high mesh resolutions. 
The time step restriction for a SIA model with an explicit or a semi-implicit discretization of the free-surface equation 
is on the form $\Delta t < C \Delta x^2$, where $\Delta t$ is the time step, and $\Delta x$ is the horizontal mesh resolution \citep{hindmarsh_payne_1996,Hindmarsh2001,Bueler_Lingle_Kallen-Brown_Covey_Bowman_2005,CHENG201729}. 
Only for extremely thin ice or steep surface gradients a linear time-step restriction occurs \citep{CHENG201729}.
A recent study showed that this quadratic behaviour carries over to hybrid models combining the SIA model with the SSA model \citep{tc-16-689-2022}. 
Resolving complex costal ice dynamics requires a fine spatial resolution in the horizontal direction, so that $\Delta x$ may be less than 1 km locally, 
and a time resolution of around $0.1$ years to $10$ years \cite{Bueler2023}. Using the SIA velocity fields, 
the simulations, however, require significantly finer time resolutions 
due to numerical instabilities, rather than physical instabilities \cite{Bueler2023}.
To alleviate the problem when considering moving ice margins, 
the SIA model was combined with a fully implicit time stepping scheme in \citep{Bueler2016}. 
This requires an implementation of a nonlinear iteration increasing 
the compuational cost and is not guaranteed to converge for bedrocks with steep gradients when the horizontal resolution is fine.

The SIA model is most commonly posed on strong form from which a closed form solution 
(explicit expressions) are obtained. 
Evaluating the closed form solution requires: (i) numerical differentiation of the ice surface position, 
(ii) numerical integration in the vertical direction, observed from the bedrock 
to the ice surface. To facilitate (ii) the mesh vertices have to be aligned over lines following the vertical direction, 
i.e., extruded meshes are needed. 
A simpler approach to implementing the SIA model in softwares based on 
Finite Element Method (FEM) is to pose the SIA model on weak form and solve the problem as a 
coupled system. This is computationally more expensive as compared to evaluating a closed form solution, 
but the weak form SIA models are still a linear problem, 
in addition allowing for fully unstructured meshes and an easier implementation using one of the FEM libraries. 
SIA models are implemented on weak form in at least two of the large scale FEM ice sheet models 
\citep{ISSM,ElmerDescrip}. 

FSSA (Free Surface Stabilization Algorithm) is an easy-to-implement, 
computationally inexpensive method for overcoming the small time-steps 
invented for mantle convection simulations \citep{KAUS201012}, 
and later introduced in the scope of ice sheet modeling 
for the nonlinear (full) Stokes problem \citep{FSSA1,FSSA2}. 
One of the requirements of the stabilization method is that the governing equations are written as a system of equations on weak form.

In this paper we consider the SIA models written as a system of equations on weak form. 
This makes it possible to add the FSSA stabilization terms. 
We discuss how the weak SIA formulation 
can best be implemented using FEM and how it can be combined with FSSA. 
We show computationally that when the SIA problem is stabilized by using the FSSA terms, 
the time step restriction is improved from quadratic to linear scaling 
in terms of the horizontal mesh size.
We further extend the weak SIA formulation to the weak linear Stokes formulation, 
which is a full Stokes model using the SIA viscosity function. 
The model does not require evaluating nonlinear iterations, 
but at the same time includes all stress components in the momentum balance. We argue 
that this improves the model robustness in terms of the numerical stability, but also 
improves the model accuracy (as compared to the weak SIA formulation) for a negligible 
increase in the computational cost. For all the enhanced SIA formulations 
we give a theoretical performance analysis estimating the operation count, 
and draw a comparison towards the operation count of the standard SIA formulation.
We focus our study on simplified two dimensional ice sheet domains: 
slab on a slope with a surface perturbation \citep{CHENG201729,Hindmarsh2001}, an idealized ice cap, and a horizontal cross-section of Greenland \citep{bedmachine}.

The paper is organized as follows. 
In Section \ref{section:governingequations} we state the different SIA and Stokes formulations 
that we consider in this paper, 
together with the free surface equation. In Section \ref{section:free_surface_equation_discretization} 
we provide information on the semi-implicit time stepping method for solving the free-surface equation. 
In Section \ref{section:sia_and_stokes_fem_discretization} we outline the spatial discretization methods 
for solving the momentum balances. 
In Section \ref{sec:fssa} we define the FSSA stabilization terms and make indications on how their addition to the SIA model 
impacts the free surface equation.
In Section \ref{section:performance_analysis} we outline a computational cost analysis of 
the considered SIA formulations. 
In Section \ref{section:experiments} we provide results to a set of numerical experiments 
assessing the time step restrictions and the error vs. runtime ratios for all the considered SIA formulations. 
In Section \ref{section:conclusion} we give our final remarks.

\section{Governing equations}
\label{section:governingequations}
In this paper we consider ice sheets that evolve in their shape as a function of time $t$. 
A simplified two dimensional ice sheet geometry is accounted for by a computational domain $\Omega=\Omega(t)$. 
One of the approaches to advance $\Omega(t)$ from time $t_k$ to time $t_{k+1}$ is to:
\begin{enumerate}
\item solve the momentum balance equations over $\Omega(t_k)$ for horizontal and vertical velocity components 
$u_1^k$ and $u_2^k$,
\item extract the ice sheet surface velocities $u_{1,s}^k$ and $u_{2,s}^k$ from $u_1^k$ and $u_2^k$ respectively,
\item solve the free surface equation using $u_{1,s}^k, u_{2,s}^k$ as data coefficients to get a new ice sheet domain $\Omega^{k+1}$.
\end{enumerate}
In this section we state the free surface equation and all the different momentum balance equations that we consider in this paper.
\subsection{Free surface equation for advancing the ice surface in time}
\label{section:governingequations:freesurfaceeq}
To compute the evolution of the ice surface function $h=h(x,t)$ in time we solve the free-surface equation:
\begin{equation}
\label{eq:free_surface_eq}
\begin{aligned}
\partial_t h &= - u_{1,s}(x,h)\, \partial_x h + u_{2,s}(x,h) + a(x,h),\quad t>0,\quad x \in \Omega^\perp,
\end{aligned}
\end{equation}
where $\Omega^\perp$ is a projected domain only taking into account the horizontal components of $\Omega$.
%
Furthermore $u_{1,s}(x,h)$ and $u_{2,s}(x,h)$ are the surface horizontal and vertical velocity functions respectively.
Term $a=a(x,h)$ is the surface mass balance in this paper set to $a(x,h) = 0$. 
We chose to work with the free surface equation rather than the thickness equation (common when using the SIA models), 
as this allows for a better flexibility in terms of 
using the free surface equation discretizations already available from the existing full Stokes model codes such as Elmer or ISSM.

The evolution of the ice surface height $h$ defines the evolution of shape of the domain $\Omega$, where $\Omega$ is 
representing the volume of an ice sheet. 
The boundary of the domain $\partial\Omega \subset \mathbb{R}$ consists of three disjoint parts:
$$\partial\Omega = \Gamma_b \cup \Gamma_s \cup \Gamma_l ,$$
where $\Gamma_b$ is the ice sheet bedrock, $\Gamma_s$ is the ice sheet free surface defined by the surface height $h$, and $\Gamma_l$ 
is the ice sheet lateral boundary. 

As is observed from \eqref{eq:free_surface_eq}, the surface $h$ is a function of the velocities $u_1$ and $u_2$. 
The velocities are computed before solving \eqref{eq:free_surface_eq}, by solving 
the momentum balance equations (SIA or Stokes) over $\Omega$. 
The coupling between $h$, $u_1$, and $u_2$ has an important impact on the time step restriction 
in the ice sheet simulations and is thus one of the main focus points of this paper. 

When solving \eqref{eq:free_surface_eq} we impose the boundary conditions as follows. We let the lateral margins of the ice sheet surface fixed, 
and use either the periodic boundary conditions: 
\begin{equation}
\label{eq:freesurface_periodic_bcs}
h(x_L,t) = h(x_R, t),
\end{equation}
where $x_L = \min_{x \in \Omega} x$ and $x_R = \max_{x \in \Omega} x$, or Dirichlet boundary conditions:
\begin{equation}
\label{eq:freesurface_dirichlet_bcs}
h(x_L, t) = 0,\quad h(x_R, t)=0.
\end{equation}
This excludes influence from nonlinearities introduced by the moving lateral boundaries, 
which is a complex problem in itself \citep{Werder,gmd-13-6425-2020,Bueler2023}. 
The margins are sometimes fixed in practice for technical reasons, see some of the models in the ISMIP6 
Antartica benchmark \citep{tc-14-3033-2020}, but clearly it is important to get the right physical 
response. 

\subsection{Strong form nonlinear full Stokes equations}
We use the full Stokes equations as a reference model for computing the velocity field 
$\mathbf{u}=(u_1,u_2)$ over an ice sheet geometry, 
when drawing the comparison towards solutions of the different SIA model formulations. 
This is reasonable as the SIA model is an approximation of the full Stokes equations. 
The full Stokes equations are:
\begin{equation}
 \begin{aligned}
    \label{eq:fullstokes}
    -\nabla \cdot \left(2 \mu_*(\mathbf{Du}) \mathbf{D}\mathbf{u} \right) + \nabla p &= \rho \mathbf{g} & \text{ on } \Omega,\\
    \nabla \cdot \mathbf{u} &= 0& \text{ on } \Omega,
\end{aligned}
\end{equation}
where $\rho>0$ is the ice density, $\mathbf{g}=(0,-9.81)$ is the gravitational acceleration, $p$ is the pressure, 
and the symmetric strain rate tensor $\mathbf{D}\mathbf{u} = \frac{1}{2} \left (\nabla \mathbf{u} + \nabla \mathbf{u}^T \right )$ is defined through four components: 
\begin{equation}
    \label{eq:Dij_strain_terms}
\begin{aligned}
    D_{11} = \partial_x u_1,\quad D_{12} = \frac{1}{2}(\partial_y u_1 + \partial_x u_2),\quad D_{21} = D_{12},\quad D_{22} =  \partial_y u_2.
\end{aligned}
\end{equation}
The viscosity function $\mu_* = \mu_*(\mathbf{Du})$ relates the strain rates to the deviatoric stress tensor $\textbf{Su}$ as 
\begin{equation}
    \label{eq:Deviatoric_stress_tensor}
\mathbf{Su}=2 \mu_*(\mathbf{D}\mathbf{u}) \mathbf{D}\mathbf{u},
\end{equation}
and is defined by:

\begin{equation}
    \label{eq:Stokesvisco}
    \mu_*(\bm D \bm u) = \frac{1}{2} A(T)^{-\frac{1}{3}} (\| \bm D\bm u\|_F^2)^{-1/3} + \varepsilon.
    \end{equation}    
%
Here $A(T)$ is constant since we consider isothermal conditions and $\varepsilon$ is the critical shear rate (a small number) 
which we define as in \citep{Hirnpowerlaw}. 

We note that the weak form of the nonlinear full Stokes problem (W-Stokes)

In all the considered test cases we impose stress free boundary conditions 
$(\mathbf{Du} - p\mathbf{I})\cdot \mathbf{n} = 0$ at the ice sheet surface $\Gamma_s$, where $\mathbf{n}$ is 
the normal vector pointing outwards of $\Gamma_s$. 
Depending on the test case we impose either the periodic or the no-slip ($\mathbf{u} = \mathbf{0}$) boundary conditions 
over the ice sheet lateral boundary $\Gamma_l$. On the ice sheet bedrock $\Gamma_b$ we impose no-slip boundary conditions 
in all test cases.

In this paper we use full Stokes equations written on weak form (abbreviated W-Stokes) defined later in the final paragraph of Section \ref{sec:governingequations:w-siastokes}. The full Stokes equations are nonlinear which leads to 
an increase computational cost when discretized and solved on a computer, 
as compared to a linear problem such as the SIA equations.






\subsection{Strong form SIA equations (SIA)}


The SIA model is derived by using an assumption that an ice sheet is thin, and that the vertical shearing forces are dominating the horizontal shearing forces. 
The stress tensor $\mathbf{Su}$ as defined in \eqref{eq:Deviatoric_stress_tensor} 
is approximated so that only the vertical shear stresses remain in place \citep{GreveBlatterBok}: 
\begin{equation}
    \label{eq:sia_stress_tensor}
\mathbf{S}=\begin{bmatrix}
 S_{11} &  S_{12}\\ 
 S_{12}& S_{22}
\end{bmatrix} \approx \begin{bmatrix}
 0 &  S_{12}\\ 
 S_{12}& 0
\end{bmatrix} = \begin{bmatrix}
 0 &  \mu \partial_y u_1\\ 
 \mu \partial_y u_1& 0
\end{bmatrix}.
\end{equation}
%
The strong form SIA equations are then:
\begin{equation}
        \label{eq:sia_strongform}
        \begin{aligned}
            - \partial_y \mu \partial_y u_1 + \partial_x p &= 0&\, \text{ on } \Omega, \\ 
             \partial_y p &= \rho g&\, \text{ on } \Omega, \\
            \partial_x u_1 + \partial_y u_2 &= 0&\, \text{ on } \Omega.
        \end{aligned}
\end{equation} 
The boundary conditions for \eqref{eq:sia_strongform} are:
\begin{equation}
\begin{aligned}
    \label{eq:sia_strong_bcs}
u_1 = 0\, \text{ on } \Gamma_b,\quad
u_2 = 0\, \text{ on } \Gamma_b,\quad
p = 0\, \text{ on } \Gamma_s,\quad
S_{12} &= 0 &\, \text{ on } \Gamma_s,
\end{aligned}
\end{equation}
As the pressure is decoupled from $u_1$ and $u_2$, we first solve the second equation of \eqref{eq:sia_strongform} 
for pressure by vertically integrating the relation from $y$ to $h(x)$:
\begin{equation}
p=\rho g (y-h),
\end{equation}
where we additionally used that $p(x,y=h)=0$.
Inserting this hydrostatic pressure in the first equation of \eqref{eq:sia_strongform} and solving for $S_{12}=\mu \partial_y u_1$ using the vertical integration gives:
\begin{equation}
S_{12}=\rho g\, \partial_x (y-h),
\end{equation}
where we also used that $S_{12}(x,y=h) = 0$.
We then compute the SIA viscosity starting at \eqref{eq:Stokesvisco}, 
using that $\|\mathbf{S}\|_F^2 = S_{12}^2$, which is due 
to the SIA assumptions defined in the scope of \eqref{eq:sia_stress_tensor}, leading to:
%
\begin{align}
    \label{eq:SIAvisco}
\mu =\frac{1}{2} A_0^{-1} (S_{12}^2 )^{-1} &= \frac{1}{2} A_0^{-1} (\rho g)^{-2} (y-h(x))^{-2} (|\partial_{x} h(x) |^{2})^{-1}  \\
& \approx  \frac{1}{2} \Big(A_0 (\rho g)^{2} (y-h(x))^{2} (|\partial_{x} h(x) |^2) + \varepsilon\big)^{-1}.
\end{align}
Here $\varepsilon$ is Hirn's regularization parameter preventing the viscosity from taking 
infinite values where the square of the surface slope $|\partial_{x} h(x) |^2 \to 0$. 
To derive \eqref{eq:SIAvisco} we also assumed isothermal conditions $A(T) = A_0 = 100$ $MPa^{-3} yr^{-1}$, where $T$ is the temperature, 
but this is generally not a limitation of the SIA model.
Using the viscosity function \eqref{eq:SIAvisco}, 
the horizontal velocity $u_1$ is given by integrating the first equation of \eqref{eq:sia_strongform} along a vertical line from $b(x)$ to $y$, and inserting that $u_1|_{b(x)}=0$. 
The vertical velocity $u_2$ is obtained by inserting the computed $u_1$ into the third equation of \eqref{eq:sia_strong_bcs}, integrating over a vertical line from $b(x)$ to $y$, 
and inserting $u_2|_{b(x)}=0$. The closed form expressions are:
\begin{equation}
    \label{eq:sia_velocities}
    \begin{aligned}
u_1 &= - \frac{1}{2}A_0(\rho g)^3 (\partial_x h)^3  \left( (y-h)^4- (b-h)^4 \right)\\
u_2 &= \frac{1}{2}A_0\, (\rho g)^3 
\Bigg(\left(\frac{1}{5}(y-h)^5 - \frac{1}{5}(b-h)^5 - (b-h)^4\,(y-b)\right)\, 3(\partial_x h)^2\, \partial_{xx} h  + ... \\
& ... - (\partial_x h)^4\, ((y-h)^4 - (b-h)^4) - 4 (\partial_x b - \partial_x h)\, (\partial_x h)^3  (b-h)^3\, \big((y-h) - (b-h) \big)\Bigg)
\end{aligned}
\end{equation}
These closed form velocity expressions are computationally inexpensive to evaluate 
as compared to solving the full nonlinear Stokes system. 
The vertical integration however requires the mesh nodes to be aligned in the vertical direction in the case of adiabatic conditions, that is, when $A$ varies with depth.

The free surface equation requires the velocities to be evaluated at the surface. Setting $y=h$ in \eqref{eq:sia_velocities} 
leads to:
\begin{equation}
    \label{eq:sia_velocities_surface}
    \begin{aligned}
u_{1,s} &= \frac{1}{2}A_0(\rho g)^3 (\partial_x h)^3  \left( (b-h)^4 \right)\\
u_{2,s} &= \frac{1}{2}A_0\, (\rho g)^3 
\left(\left(- \frac{3}{5}(b-h)^5 + 3(b-h)^5 \right)\, (\partial_x h)^2\, \partial_{xx} h  + \left(4(b-h)^4\, \partial_x b\, (\partial_x h)^2 - 3(b-h)^4(\partial_x h)^3 \right)\, \partial_x h \right)
\end{aligned}
\end{equation}

After inserting the SIA velocities from \ref{eq:sia_velocities} to the free surface equation 
\ref{eq:free_surface_eq}, 
we write the free-surface equation problem as a nonlinear advection-diffusion PDE:
\begin{equation}
\begin{aligned}
    \label{eq:free_surface_eq_advection_diffusion}
    \partial_t h &= - u_{1,s}\,\partial_x h + u_{2,s} \\
    &= C_1(h)\, \partial_x h + C_2(h)\, \partial_{xx} h,
\end{aligned}
\end{equation}
where:
\begin{equation}
    \label{eq:free_surface_eq_advection_diffusion_coefficients}
    \begin{aligned}
        C_1(h) &=  \frac{1}{2}A_0(\rho g)^3 \left(  (b-h)^4\, (\partial_x h)^2 + 4(b-h)^4\, \partial_x b\, (\partial_x h)^2 - 3(b-h)^4(\partial_x h)^3 \right)   \\
        C_2(h) &=  \frac{1}{2}A_0\, (\rho g)^3 \left(- \frac{3}{5}(b-h)^5\, (\partial_x h)^2 + 3(b-h)^5\, (\partial_x h)^2 \right)\,  .\\
    \end{aligned}
    \end{equation}
The time step restriction has a quadratic scaling in terms of the mesh size 
due to the second derivative term (diffusive term) $\partial_{xx} h$ in \eqref{eq:free_surface_eq_advection_diffusion}. 
The standard way to theoretically assess the timestep restriction is to linearize 
\eqref{eq:free_surface_eq_advection_diffusion} with respect to $h$ and then perform 
a von Neumann (Fourier) analysis. 
This was done in e.g. \citep{CHENG201729} 
for the thickness equation in the case of a perturbed slab on a slope. 
We repeat this exercise for the free surface equation in the appendix and will revisit it for a new SIA formulation where the FSSA stabilization of  \citep{FSSA1,FSSA2} is added. 

\subsection{Weak form SIA equations (W-SIA)}
The easiest approach to implementing the SIA equations within an existing FEM code -- 
such as Elmer, ISSM, or FEniCS -- is to write 
\eqref{eq:sia_strongform} on weak form and discretize the weak form using the standard finite element methods. 
This also allows for fully unstructured meshes, which can be of higher quality on certain geometries, and is sometimes techincally easier to construct. 
The weak SIA formulation is obtained by multiplying each equation of 
\eqref{eq:sia_strongform} using piecewise continuous test functions 
$v_1=v_1(x,y),v_2=v_2(x,y)$, $q=q(x,y)$, respectively, and integrating over $\Omega$. 
The first term of the first equation is additionally integrated by parts. In the end the weak SIA formulation 
is:
\begin{equation}
    \label{eq:sia_weakform_less}
    \begin{aligned}
        \int_{\Omega} \mu \partial_y u_1\, \partial_y v_1\, d\Omega - \int_\Omega \partial_x p\, v_1\, d\Omega &= 0, \\ 
         - \int_{\Omega} \partial_y p\, v_2\, d\Omega &= \int_\Omega \rho g\, v_2\, d\Omega, \\
        \int_\Omega \left(\partial_x u_1 + \partial_y u_2\right)\, q\, d\Omega &= 0,
    \end{aligned}
\end{equation}
where $\mu$ is the SIA viscosity defined in \eqref{eq:SIAvisco}. 
In this paper we abbreviate the weak SIA formulation as W-SIA.
Solving W-SIA on a computer is cheaper as compared to solving the 
full nonlinear Stokes system \eqref{eq:fullstokes}, since W-SIA is a linear problem 
due to that the viscosity is not a function of the computed velocity. 
W-SIA is possible to solve in terms of three subsequent matrix systems: 
first for pressure, secondly for $u_1$, and lastly for $u_2$. This only holds as long as W-SIA is not further stabilized using the 
additional stabilization terms that couple the velocity functions.
A disadvantage when solving W-SIA is that 
the many stress components are not present in \eqref{eq:sia_weakform_less}. 
This implies that the full stress term $\mathbf{Su}$ is not guaranteed to have an upper bounded 
when the problem is solved on a computer and the mesh size goes to $0$. A consequence are potentially 
sharp velocity gradients that deteriorate the numerical stability as well as the solution accuracy.

\subsection{Weak form linear Stokes equations employing the SIA viscosity function (W-SIAStokes)}
\label{sec:governingequations:w-siastokes}
We add the missing stress components back to \eqref{eq:sia_weakform_less} resulting in the following 
weak formulation:
\begin{equation}
    \label{eq:sia_weakform_all}
    \begin{aligned}
        \int_{\Omega} 2\mu \partial_x u_1\, \partial_x v_1\, d\Omega +  \int_{\Omega} \mu (\partial_y u_1 + \partial_x u_2)\, \partial_y v_1\, d\Omega - \int_\Omega \partial_x p\, v_1\, d\Omega &= 0, \\ 
        \int_{\Omega} \mu (\partial_y u_1 + \partial_x u_2)\, \partial_x v_2\, d\Omega +  \int_{\Omega} 2\mu \partial_y u_2\, \partial_y v_2\, d\Omega - \int_{\Omega} \partial_y p\, v_2\, d\Omega &= \int_\Omega \rho g\, v_2\, d\Omega, \\
        \int_\Omega \left(\partial_x u_1 + \partial_y u_2\right)\, q\, d\Omega &= 0.
    \end{aligned}
\end{equation}
We abbreviate \ref{eq:sia_weakform_all} as W-SIAStokes, as that formulation combines the full Stokes problem and the SIA viscosity function. 
 W-SIAStokes is a linear problem, computationally less expensive to solve as compared to the 
 full nonlinear Stokes problem. 
However, the system can no longer be solved as three seperate matrix systems as is the case in the unstabilized W-SIA formulation \eqref{eq:sia_weakform_less}. 
An advantage of W-SIAStokes over W-SIA is a guaranteed bound over the full stress term $\mathbf{Su}$, which 
improves the numerical stability properties as the mesh size goes to $0$.

We note that the nonlinear full Stokes problem \eqref{eq:fullstokes} 
but written on weak form (W-Stokes) takes exactly the same form as W-SIAStokes \eqref{eq:sia_weakform_all}, 
where we use the (full) viscosity function \eqref{eq:Stokesvisco} instead of the SIA viscosity \eqref{eq:SIAvisco}.
\section{Finite difference discretization of the free surface equation}
\label{section:free_surface_equation_discretization}
We first denote that $h=h(x,t)$, and discretize the free surface equation \eqref{eq:free_surface_eq} in time 
using the first order semi-implicit Euler method. This results in:
\begin{equation}
	\label{eq:discreteicesurface_time}
		\frac{h^{k+1} - h^k}{\Delta t } = - u_{1,s}^k  \partial_x h^{k+1} + u_{2,s}^k,\qquad k=1,2,3,..
\end{equation}
where $h^k, h^{k+1}$ are $h(x,t_k)$, $h(x,t_{k+1})$ respectively, and where $u_{1,s}^k$, $u_{2,s}^k$ are 
the surface velocities $u_1(x^k,y_s^k)$, $u_2(x^k,y_s^k)$ extracted from the bulk velocity functions defined over an ice sheet domain $\Omega^k$ at $t_k$. 
We note that $x \in \Omega^\perp$, where this domain is defined in the scope of Section \ref{section:governingequations:freesurfaceeq}.
Now we discretize the spatial derivatives in \eqref{eq:discreteicesurface_time} by means of 
the second-order accurate centered finite difference stencil weights, resulting in the following system of equations:
\begin{equation}
    \label{eq:free_surface_time_discretization}
		\frac{\bm h^{k+1} - \bm h^k}{\Delta t } = - \text{diag}(\bm u_{1,s}^k)  \bm D_x \bm h^{k+1} + \bm u_{2,s}^k,\qquad k=1,2,3,..
\end{equation}
where $h^{k+1}_i = h(x_i,t_{k+1})$, 
$h^{k}_i = h(x_i,t_{k})$, 
$(u_1^{k+1})_i = u_1(x_i, y_s)$, 
$(u_2^{k+1})_i = u_2(x_i, y_s)$, $i=1,..,N$. The components of the matrix $\bm D_x$ are defined by the second-order accurate finite difference stencil weights 
that discretize the first order derivative operator.
The final time stepping iteration scheme is:
\begin{equation}
    \label{eq:free_surface_time_discretization}
		\bm h^{k+1} = (\bm I + \Delta t\, \text{diag}(\bm u_{1,s}^k) \bm D_x)^{-1}\,  (\bm h^k + \Delta t\,\bm u_{2,s}^k),\qquad k=1,2,3,..
\end{equation}
We impose the boundary conditions as described within the scope of 
Section \ref{section:governingequations:freesurfaceeq} 
by reducing the system of equations \eqref{eq:free_surface_time_discretization} 
in the unknowns related to the Dirichlet conditions, or by 
transforming the $\bm D_x$ matrix into a circulant matrix in the case when we use the periodic boundary conditions.

To instead use a fully implicit scheme would require access to $u_1^{k+1}$, $u_2^{k+1}$, 
but this is computationally expensive as the velocity functions and the surface position 
are coupled. The surface $h$ depends on $u_1^k$, $u_2^k$ due to 
\eqref{eq:discreteicesurface_time}, while the velocities depend on the surface that determines the 
shape of the computational domain $\Omega$ on which we solve the momentum balance equations. 
As a consequence, computing $u_1^{k+1}$, $u_2^{k+1}$ 
requires an expensive nonlinear iteration as demonstrated in the SIA model case in \citep{Bueler2016}.

\section{Finite element discretization of the SIA / Stokes models}
\label{section:sia_and_stokes_fem_discretization}
Throughout the paper we use FEM not only to solve partial differential equations on weak form, 
but also to evaluate the surface gradient functions. 
The meshes we use are extruded. 
To create a two-dimensional ice sheet mesh we first generate a rectangular 
mesh with dimensions $[x_{\text{min}}, x_{\text{max}}] \times [0, 1]$, 
where $x_{\text{min}}$ and $x_{\text{max}}$ are the 
minimum and the maximum horizontal coordinates of the ice sheet geometry. 
Then we transform the vertical mesh coordinates using an ice sheet initial surface function.

When evaluating the SIA velocities using the closed form expression 
\eqref{eq:sia_velocities} we employ FEM to evaluate $\partial_x h$, 
the gradient of the ice sheet surface. 
We first interpolate $h$ into a piecewise linear finite element space. 
After that we compute the gradient $\partial_x h|_{K_i}$, $i=1,..,N$ over each mesh element $K_i$. 
As the gradient of the picewise linear function across 
the element interfaces is discontinuous (not well defined), 
we $L_2$ project the computed gradients back into a piecewise linear finite element space. By that we 
compute a continuous (well defined) surface gradient.

When solving the nonlinear Stokes problem \eqref{eq:fullstokes} on weak form we use Taylor-Hood elements (P2P1) in order to fulfill the 
inf-sup condition \citep{Babuska1973,Brezzi1974}, 
that is, picewise quadratic polynomials for approximating the velocity functions, 
and piecewise linear polynomials for approximating the pressure function. This is a requirement to make the finite element discretization 
numerically stable.

When solving W-SIAStokes \eqref{eq:sia_weakform_all} we use the same type 
of elements as in the nonlinear Stokes problem case, for the very same reasons related to numerical stability.

When solving W-SIA the inf-sup condition does not need to be fulfilled, and we can therefore use piecewise linear finite elements (P1P1) for approximating the velocity functions, as well as the pressure function. 
This is an advantage as the amount of unknowns when using P1P1 elements is three times smaller per an unknown function, as compared to 
when using P2P1 elements.

\section{Free Surface Stabilization Algorithm (FSSA) for the SIA / Stokes models}
\label{sec:fssa}

In \citep{FSSA1,FSSA2} the authors introduced FSSA for the full nonlinear Stokes model 
to mimic an associated implicit solver advancing the ice surface from time $t_k$ to time $t_{k+1}$, $k=1,..,N$, where $t_{k+1} > t_k$. 
This is done by predicting the gravitational force on weak form at $t_{k+1}$ by adding an extra surface force term: 
\begin{equation}\label{eq:fssa}
\int_{\Omega^{k+1}} \rho g v_2 \approx \int_{\Omega^{k}} \rho g v_2 + \theta \Delta t \int_{\partial \Omega^k} (\mathbf{u \cdot n})\, \rho g v_2\, ds.
\end{equation}
Here $\theta \in [0,1]$ is a user-defined constant parameter. The relation \eqref{eq:fssa} was derived from a finite difference 
discretization of the Reynolds transport theorem (the multi-dimensional Leibniz rule). 
The FSSA thus relies on the assumption that the flow is predominantly gravity driven, 
so that computing the gravitational force at $t_{k+1}$ leads to a good approximation of the ice flow at $t_{k+1}$. 
Hence an implicit discretization \citep{FSSA1,FSSA2}. From a physics standpoint, 
the FSSA term is an extra surface pressure term acting as a damping term -- when the ice is rising, the FSSA term acts as an extra surface pressure, and when the ice is sinking, it reduces the pressure. 
FSSA was originally introduced by \citep{KAUS201012} for mantle convection simulations.

In this paper we add the FSSA stabilizion term to to the vertical momentum balance equation. 
In the W-SIAStokes case \eqref{eq:sia_weakform_all} (and similar in the W-Stokes case) the vertical momentum balance equation becomes:
\begin{equation}
  \int_{\Omega} \mu (\partial_y u_1 + \partial_x u_2)\, \partial_x v_2\, d\Omega +  \int_{\Omega} 2\mu \partial_y u_2\, \partial_y v_2\, d\Omega - \int_{\Omega} \partial_y p\, v_2\, d\Omega = \int_\Omega \rho g\, v_2\, d\Omega + \theta \Delta t \int_{\partial\Omega} (\mathbf{u \cdot n})\, \rho g v_2\, ds.
\end{equation}
In the W-SIA case \eqref{eq:sia_weakform_less} the vertical momentum balance, after adding the FSSA stabilization term, becomes:
\begin{equation}
        - \int_{\Omega} \partial_y p\, v_2\, d\Omega = \int_\Omega \rho g\, v_2\, d\Omega + \theta \Delta t \int_{\partial\Omega} (\mathbf{u \cdot n})\, \rho g v_2\, ds.
\end{equation}
In this case it is the added FSSA term that couples the pressure to surface velocities $\mathbf{u}_s$. 
Without the FSSA term the pressure is decoupled from the velocity, reducing the computational cost of the solution procedure. 
The coupling is however essential for numerical stability reasons. 
The FSSA term renders the pressure implicit, and by that allows for long time steps when solving the free surface equation. 
To show that, we start at the pressure from the strong SIA formulation \eqref{eq:sia_strongform}, and evaluate 
it at $t_{k+1}$:
\begin{align}
         p^{k+1} &= \rho g (h^{k+1}-y)\\
         \label{eq:FSSAmakesPimplicit_step1}
\end{align}
Adding and subtracting $p^k = \rho gh^k$, and then reordering the terms, gives:
\begin{equation}
    p^{k+1} = \rho g  (h^k -y)+ \rho g\, (h^{k+1} -h^k).
    \label{eq:FSSAmakesPimplicit_step2}
  \end{equation}
Now we observe that the explicit Euler discretization of the free surface equation \eqref{eq:free_surface_eq} is 
$h^{k+1} - h^k = \Delta t \rho g(- u_{1,s}\partial_x h + u_{2,s})$. 
We insert the relation to \eqref{eq:FSSAmakesPimplicit_step1}, and obtain:
\begin{equation}
  \label{eq:FSSAmakesPimplicit_step2}
  \begin{aligned}
  p^{k+1} &= \rho g  (h^k -y)+ \Delta t \rho g(- u_{1,s}\partial_x h + u_{2,s}) \\
  &= p^k + \Delta t \rho g(- u_{1,s}\partial_x h + u_{2,s}).
  \end{aligned}
\end{equation}  
Pressure at $t_{k+1}$ is thus the pressure at $t_k$ plus the FSSA term (on strong form). 
The same relation as \eqref{eq:FSSAmakesPimplicit_step2} 
can also be derived by applying the Leibniz rule to the vertical integration of the second equation of \eqref{eq:sia_strongform}. 
We show in appendix \ref{sec:appendix:von_neumann} for the strong form SIA \eqref{eq:sia_strongform} 
that assuming an implicit pressure (i.e. pressures $p^{k+1}$ evaluated at $t_{k+1}$), then this is sufficient to alleviate 
the quadratic time step restriction $\Delta t < C \Delta x^2$ when solving the free surface equation. To derive this result we 
extend the Von Neumann type analysis for a slab on a slope test case from \citep{CHENG201729} to our case. 
In \citep{CHENG201729} the authors show that, assuming thick ice with low surface inclination, 
the quadratic dependence on $\Delta x$ is:
\begin{equation}
\Delta t<\frac{3}{5}A_0 |\rho g|^3  C_\alpha^2  \bar{H}^5 (\Delta x)^2,
\end{equation}
where $C_\alpha$ is the average surface slope and $H$ the ice thickness.
The result stems from that the vertical velocity $u_2$ contains a second derivative of the surface $\partial_{xx} h$. 
Furthermore, the second derivative of the surface origins from that the vertical velocity is a function of 
$\partial_x u_1$ which is in turn 
a function of the horizontal pressure derivative $\partial_x p= \rho g\, \partial_x h^k$. 
Following the derivation of the compact form free surface equation \eqref{eq:free_surface_eq_advection_diffusion} 
with the coefficients \eqref{eq:free_surface_eq_advection_diffusion_coefficients}, we now 
write the time discretized free surface equation when assuming that the velocities are derived from $p^{k+1}$, 
but that the intermediate integration steps involve $h^k$. We have:

\begin{equation}
      \label{eq:fssa_free_surface_time_discretized}
      h^{k+1} = h^k + C_0\, \partial_x h^{k} + C_1\, \partial_x h^{k+1} + C_2\, \partial_{xx} h^{k+1}
  \end{equation}
  where:

  \begin{equation}
    \label{eq:fssa_free_surface_time_discretized_coefficients}
    \begin{aligned}
      C_0 &=  \frac{1}{2}A_0(\rho g)^3 \left( - 3(b-h^k)^4(\partial_x h^{k+1})^3 \right)   \\
      C_1 &=  \frac{1}{2}A_0(\rho g)^3 \left(  (b-h^k)^4\, (\partial_x h^{k+1})^3  + 4(b-h^k)^4\, \partial_x b\, \partial_x h^{k+1} \right)   \\
      C_2 &=  \frac{1}{2}A_0\, (\rho g)^3 \left(- \frac{3}{5}(b-h^k)^5\, (\partial_x h^{k+1})^2 + 3(b-h^k)^5\, (\partial_x h^{k+1})^2 \right)\,  .\\
    \end{aligned}
    \end{equation}
We now have an implicit treatment of the leading diffusive $\partial_{xx} h$ term in \eqref{eq:fssa_free_surface_time_discretized}
 which leads to a linear time step constraint:
\begin{equation}
  \Delta t \leq \Big(\frac{3}{2}A_0 |\rho g|^3 C_\alpha^3  \bar{H}^4 \Big)^{-1}\, \Delta x,
  \end{equation}  
where $C_\alpha$ is the average surface slope. 
We derive the above time step restriction for 
the slab on a slope with a perturbed ice surface case in Appendix \ref{sec:appendix:von_neumann} 
and validate the estimate numerically for W-SIA and W-SIAStokes in the numerical experiments section.




\section{Computational cost estimation for the different SIA / Stokes formulations}
\label{section:performance_analysis}
The computational work when solving momentum balance models is highly dependent on both software and hardware. 
However, it is still possible to make estimates of the computational cost, see for instance 
a type of performance analysis approach of \citep{Bueler2023}. 
In this section we make rough estimates of the computational cost for ice sheet simulations, 
where the velocity functions are computed using SIA \eqref{eq:sia_strongform}, W-SIA \eqref{eq:sia_weakform_less}, W-SIAStokes \eqref{eq:sia_weakform_all}, W-Stokes \eqref{eq:Stokesvisco}, 
and the ice surface is advanced from time $t_k$ to time $t_{k+1}$, $k=1,..,N_{\Delta t}$, 
using the discretized free surface equation \eqref{eq:free_surface_time_discretization}. 
We write the approximate computational cost on the form:
$$\text{Computational cost} = C(d,\alpha)\, C_S\, m^{1 + \gamma/(d-1) + \alpha },$$
where:
\begin{itemize}
   \item $m$ is the number of mesh vertices (nodes) in the horizontal direction,
   \item $\alpha \in [0, 2]$ denotes the choice of a linear solver ($\alpha = 2$ dense direct solver, $\alpha = 1$ sparse direct solver, 
   $\alpha = 0.05$ algebraic multigrid solver \citep{Bueler2023}). For pure SIA no linear solver is needed and $\alpha=0$. 
   \item $\gamma$ is the scaling exponent in the simulation time step restriction $\Delta t  \leq C_t\, \Delta x^\gamma$, where $\Delta x$ is the horizontal internodal distance.
   \item $C_{S}$ is a constant specific to the computational cost of the nonlinear Stokes problem which 
   involves the nonlinear iteration count and the choice of hardware,
   \item $C(d,\alpha)$ is a constant depending on $\alpha$ and $d$, where $d$ is the dimension count of the considered ice sheet geometry.
  \end{itemize}
Following \citep{Bueler2023} we have that W-Stokes requires $C_{S}\, m^{1+\alpha}$ floating point operations until convergence 
per one time step.
%
%

In the W-SIAStokes case the computational cost is the same as in W-Stokes, 
but divided by the nonlinear iteration count $N_{\text{iter}}$. The cost is $\frac{1}{N_{\text{iter}}}\, C_{S}\, m^{1+\alpha}$. 
This is due to that W-SIAStokes only differs from W-SIA in the choice of the viscosity function 
(linear) \eqref{eq:SIAvisco} and thus requires one iteration to be solved. When making the estimate 
we also assumed that the different viscosity function preserves the preconditioning 
quality. 

When W-SIA \eqref{eq:sia_weakform_less} is not FSSA stabilized then the three ($d+1$ equations in general) equations 
are solved one by one, and so in this case $m \to \frac{m}{d+1}$. This gives the estimate 
$\frac{d+1}{(d+1)^{1+\alpha}}\, \frac{1}{N_{\text{iter}}}\, C_{S}\, m^{1+\alpha}$. When W-SIA is FSSA stabilized then the decoupled 
solution procedure is not possible to perform anymore and the computational cost is the same as in the W-SIAStokes case, that is, 
$\frac{1}{N_{\text{iter}}}\, C_{S}\, m^{1+\alpha}$.

Computing the SIA velocities by means of the closed-form expressions \eqref{eq:sia_strongform} requires $C_{\text{SIA}}\, m$ 
floating point operations.

The computational cost for advancing the ice surface from the initial state to time $t=T$ 
is proportional to the number of time steps $N_{\Delta t}$ we take along the way. 
The number of time steps itself is given by $N_{\Delta t} = \frac{T}{\Delta t} \sim \frac{1}{\Delta t}$. 
For a time step restriction on the form $\Delta t  \leq C_t\, \Delta x^\gamma$, we have that 
$N_{\Delta t} \sim \frac{1}{\Delta x^\gamma} \sim m^{\gamma/(d-1)}$. 

We now combine the computational cost estimates for obtaining the velocity functions, with the computational cost estimate for advancing the ice surface 
in time. The estimates for all the considered formulations are gathered in the table below.
%
\\
\begin{center}
{\tabulinesep=0.6mm
\begin{tabu}{ l || c  }
  \textbf{Model} &  \textbf{Computational cost estimate} \\ \hline \hline 
  \textbf{W-Stokes} &  $C_S\,  m^{1+\gamma/(d-1)+\alpha}$   \\ \hline
  \textbf{W-Stokes-FSSA} &  $C_S\,  m^{1+\gamma/(d-1)+\alpha}$   \\ \hline
  \textbf{W-SIAStokes} & $\frac{1}{N_{\text{iter}}}\, C_{S}\, m^{1+\gamma/(d-1)+\alpha}$ \\ \hline
  \textbf{W-SIAStokes-FSSA} & $\frac{1}{N_{\text{iter}}}\, C_{S}\, m^{1+\gamma/(d-1)+\alpha}$ \\ \hline
  \textbf{W-SIA} & $\frac{d+1}{(d+1)^{1+\alpha}}\, \frac{1}{N_{\text{iter}}}\, C_{S}\, m^{1+\gamma/(d-1)+\alpha}$  \\ \hline
  \textbf{W-SIA-FSSA} & $\frac{1}{N_{\text{iter}}}\, C_{S}\, m^{1+\gamma/(d-1)+\alpha}$ \\ \hline
  \textbf{SIA} & $C_{\text{SIA}}\, m^{1+\gamma/(d-1)}$  \\
\end{tabu}
}
\end{center}
\, \\
All the parameters to evaluate the computational costs in the above table are known, 
except for the time step restriction exponent $\gamma$
in the case of some SIA formulations. We numerically compute the exponents $\gamma$ in Section \ref{section:experiments}, 
and then compare the computational cost across the different formulations.
\section{Numerical study}
\label{section:experiments}
In this section we solve the SIA problem posed on: the strong form \eqref{eq:sia_strongform}, 
the reduced weak form \eqref{eq:sia_weakform_less}, and the 
full weak form \eqref{eq:sia_weakform_all}. 
We numerically compute the largest stable time step size $\Delta t$ when the free surface of an ice sheet is 
advected in time as described in Section \ref{sec:computations_largest_timestep}. We 
find the dependence between $\Delta t$ and the horizontal mesh size $\Delta x$(the CFL condition), compare the errors of the different SIA solutions 
to the nonlinear Stokes solution and relate them to runtimes. We do this for three different geometries. 
The experiments are performed by using the FEniCS 2019 library \citep{UFL,Alnaes2014,Alnaes2015} on a laptop with the AMD Ryzen 7 PRO 6850U processor and  16 GB RAM.

\subsection{An algorithm for computing the largest feasible time step when solving the free surface equation}
\label{sec:computations_largest_timestep}
In this section we provide the criterion under we use to numerically compute the largest feasible time step $\Delta t_*$ 
when updating the ice sheet surface in time using the free surface equation
\eqref{eq:free_surface_eq} over domain $\Omega^\perp$ as defined in the scope of Section \ref{section:governingequations:freesurfaceeq}. 
Courant-Friedrich-Levy (CFL) condition limits $\Delta t_*$ in terms of the mesh size $\Delta x$:
\begin{equation}
    \label{eq:methods:time_step_restriction}
    \Delta t_* = C\, \min_{x_i \in \Omega^\perp} (\Delta x)^p,\, p>0,
\end{equation}
where $C>0$ is the CFL number depending on the type of the discretization of \eqref{eq:free_surface_eq} and the data in \eqref{eq:free_surface_eq}, 
and the mesh size is:
\begin{equation}
\Delta x = \min_{(x_i \neq x_j) \in \Omega^\perp_h} \|x_i - x_j\|_2,
\end{equation}
We are especially interested in the exponent $p$ from \eqref{eq:methods:time_step_restriction}, where the severity of the time step restriction 
increases with an increased $p$, whereas $p=0$ implies no dependence of $\Delta t_*$ on $\Delta x$.
To compute $\Delta t_*$ numerically we use the stability criterion:
\begin{equation}
    \label{eq:experiments:dt_energy_stability_criterion}
    \int_{\Omega^\perp} (h(x,t_{k+1}))^2\, d\Omega^\perp - \int_{\Omega^\perp} (h(x,t_{k}))^2\, d\Omega^\perp \leq 0,\quad k=1,2,3,..
\end{equation}
which has to hold for each time $t_k < T$, where $k=1,2,3,...$. The above stability criterion measures the difference in the 
energy of the free surface function across two consecutive time samples. This is motivated by the fact that in 
the Von Neumann analysis, $h=h(x,t)$ is decomposed into Fourier modes 
$\delta_j^k$, $j=1,..,N$, $k=1,2,3...$, where $j$ is the number of the mode, and $k$ is related to the time sample $t_k$. 
The final statement of the Von Neumann analysis is that $\exists \Delta t>0$ such that:
$$|\delta_j^{k+1}| \leq |\delta_j^{k}|.$$ 
Thus, $\sum_{j=1}^N |\delta_j^{k+1}| \leq \sum_{j=1}^N |\delta_j^{k}|$, 
and using Parseval's identity $\sum_{j=1}^N |\delta_j^{k}| = \int_{\Omega^\perp} h(x,t)^2\, d\Omega^\perp$ 
on each side of the inequality, and then moving the right-hand-side term to the left-hand-side we obtain 
\eqref{eq:experiments:dt_energy_stability_criterion}.

We pose the computation of $\Delta t_*$ as an optimization problem:
\begin{equation}
    \label{eq:methods:time_step_restriction_optimization}
    \begin{aligned}
\text{Find }\max \Delta t_* & \\
\text{subject to } & \int_{\Omega^\perp} (h(x,t_{k+1}))^2\, d\Omega^\perp - \int_{\Omega^\perp} (h(x,t_k))^2\, d\Omega^\perp < 0,\quad k=1,2,...,N_T,
    \end{aligned}
\end{equation}
where $\int_{\Omega^\perp} (h(x,t_{k+1}))^2\, d\Omega^\perp$ and $\int_{\Omega^\perp} (h(x,t_k))^2\, d\Omega^\perp$ 
are the free-surface energies evaluated in two consecutive time steps, and is $N_T$ is a number of time steps to perform the simulations in time $t \in (0, T]$. 
To understand how $\Delta t_*$ depends on $\Delta x$, 
we discretize \eqref{eq:free_surface_eq} using different mesh sizes $(\Delta x)_j$, $j=1,2,..$, and then for each $(\Delta x)_j$ 
compute $(\Delta t_*)_j$ by solving one optimization problem \eqref{eq:methods:time_step_restriction_optimization}. 
We solve the optimization problem using the bisection method. 
Once all the data paris $((\Delta x_*)_j, \Delta t_*)_j)$, $j=1,2,...$, are known, 
we approximate the exponent $p$ in \eqref{eq:methods:time_step_restriction} by 
fitting a linear function to the transformed data pairs $((\log_{10} \Delta x_*)_j, (\log_{10} \Delta t_*)_j)$, $j=1,2,...$, 
where $p$ takes the value of the slope of the fitted linear function.

Note that this is a fairly computationally expensive procedure, which is the reason why we restrict ourselves to two dimensions in this study. 

\subsection{Slab on a slope with perturbed surface}
\label{sec:experiments:slabslope}
\subsubsection{Configuration}
First we run experiments for a slab on a slope with a perturbed surface, as this is the setting of the von Neumann analysis in \eqref{sec:appendix:von_neumann}. The ice sheet is a two-dimensional slab, $L=80\cdot 10^3$ meters long, and $H=1 \cdot 10^3$ meters thick, inclined with $\alpha=0.75$ 
degrees measured in the clockwise direction, with the initial surface:
$$h(x,0) = H + e^{-5 \cdot 10^{-8} (x - \frac{L}{2})^2}.$$

For the free-surface equation \eqref{eq:free_surface_eq} we impose periodic boundary conditions \eqref{eq:freesurface_periodic_bcs} as this is required for the von Neumann analysis. We impose no-slip boundary conditions on $\Gamma_b$, stress-free boundary conditions on $\Gamma_s$, and 
periodic boundary conditions on $\Gamma_l$.

\begin{figure}[h!]
    \centering
\begin{tabular}{c}
    \textbf{Slab on a slope surface case} \\
\textbf{Surface elevation propagation in time (reference solution)} \\    
\includegraphics[width=0.5\linewidth]{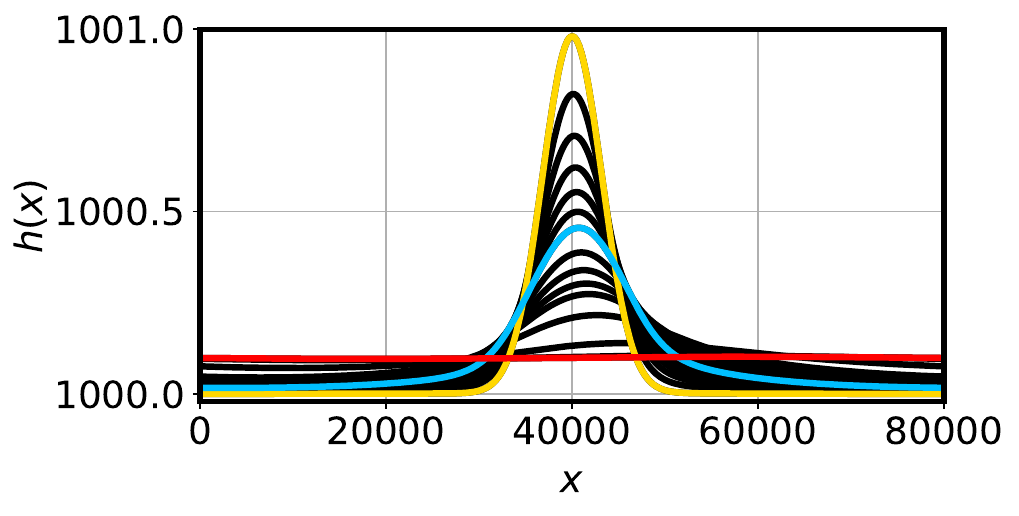}
\hspace{-0.4cm}\raisebox{2.25cm}{\includegraphics[width=0.17\linewidth]{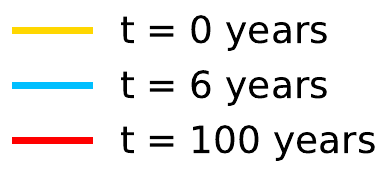}}
\end{tabular}
\caption{Propagation of the surface elevation function in time when computed as a solution to the nonlinear Stokes problem (reference), where the simulation time step is $\Delta t = 0.1$ years. 
Horizontal mesh size and vertical mesh size are $\Delta x = 250$ meters and $\Delta y = 90$ meters respectively.}
\label{fig:slabslope:solution_reference_pstokes}
\end{figure}
The relation between the number of discretiation points in the horizontal direction $n_x=m$ and $\Delta x$ that we use to perform the experiments is given in the table below. 
Note that $n_y=11$ (corresponding to $\Delta y = 90$ meters) is fixed for all experiments and the FSSA scaling parameter is fixed at $\theta = 1$, unless stated otherwise.
\\
\vspace{0.2cm}
\begin{tabular}{c|cccccccccc}
    \multicolumn{11}{c}{\textbf{Slab on a slope surface case}} \vspace{0cm} \\
    \multicolumn{11}{c}{\textbf{Number of elements and the mesh size}} \vspace{0.2cm} \\
$\mathbf{n_x}$ & 20 & 30 & 40 & 50 & 60 & 70 & 80 & 120 & 160 & 200  \\ \hline
$\mathbf{\Delta x}$ & 4000.0 & 2666.7 & 2000.0 & 1600.0 & 1333.3 & 1142.9 & 1000.0 & 666.7 & 500.0 & 400.0 \vspace{0.3cm} \\
$\mathbf{n_x}$ & 240 & 280 & 320 & 480 & 640 & 1000 & 1280 & 1560 & 1800 \\ \hline
$\mathbf{\Delta x}$ & 333.3 & 285.7 & 250.0 & 166.7 & 125.0 & 80.0 & 62.5 & 51.3 & 44.4
\end{tabular}
\vspace{0.3cm}

\begin{figure}[h!!]
    \centering
\begin{tabular}{c}
\textbf{Slab on a slope surface case} \vspace{0.2cm} \\
\textbf{Surface elevations at $\mathbf{t=6}$ years (SIA solutions)} \\    
\includegraphics[width=0.5\linewidth]{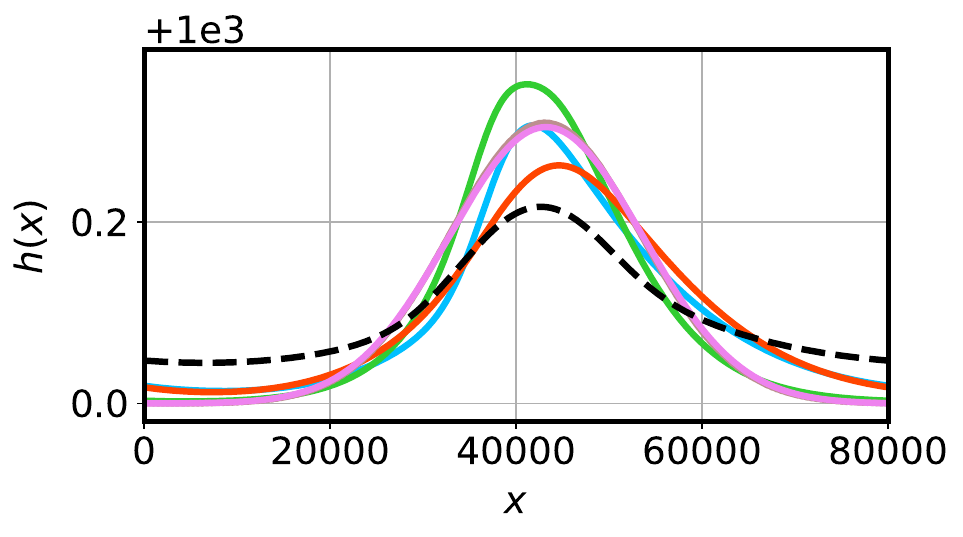}
\hspace{-0.4cm}\raisebox{1.25cm}{\includegraphics[width=0.21\linewidth]{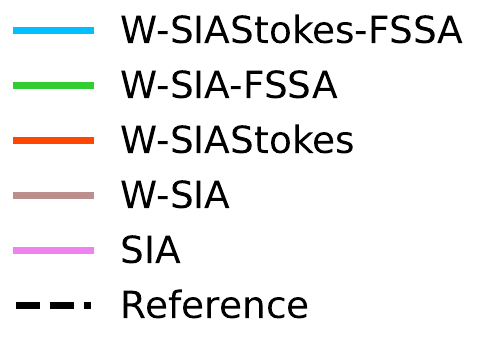}}
\end{tabular}
\caption{
    Surface elevations at simulation time $t=6$ years for different SIA problem formulations and the reference formulation (nonlinear Stokes problem). Horizontal and vertical mesh sizes are 
$\Delta x = 250$ meters and $\Delta y = 90$ meters respectively. The largest feasible time step $\Delta t_*$ for the given $\Delta x$ and $\Delta y$ 
is chosen for each formulation: 
$\Delta t_* = 6$ years, $\Delta t_* = 12$ years, $\Delta t_* = 1.8$ years, $\Delta t_* = 0.04$ years, $\Delta t_* = 0.008$ years, $\Delta t = 0.1$ years 
for (W-SIAStokes-FSSA), (W-SIA-FSSA), (W-SIAStokes), (W-SIA), (SIA), (Reference)  respectively.
}
\label{fig:slabslope:solutions_allforms}
\end{figure}

\subsubsection{Accuracy}
In Figure \ref{fig:slabslope:solution_reference_pstokes} we show how the perturbed surface evolves in time, 
where the final simulation time is $t=100$ years. 
Here the solution is computed using the nonlinear Stokes problem which we consider a reference, with a small time step $\Delta t=0.1$ years, 
the horizontal mesh size $\Delta x=250$ meters and the vertical mesh size $\Delta y = 90$ meters. 
In Figure \ref{fig:slabslope:solutions_allforms} we plot the solutions of the different SIA problem formulations to the reference solution, 
where all the solutions are evaluated at $t=6$ years. We observe that all the solutions are close to the reference solution. The solution to 
W-SIAStokes appears overall closest to the reference. The addition of the FSSA stabilization term increases the error, but 
not significantly.

\begin{figure}[h!]
    \centering
\begin{tabular}{ccc}
    \multicolumn{3}{c}{\hspace{-2cm}\textbf{Slab on a slope surface case}} \vspace{0.2cm} \\
    \hspace{1cm}\textbf{$\mathbf{\Delta t}$ vs. $\mathbf{\Delta x}$} & \hspace{1.2cm}\textbf{Error vs. Runtime} & \\    
     \includegraphics[width=0.26\linewidth]{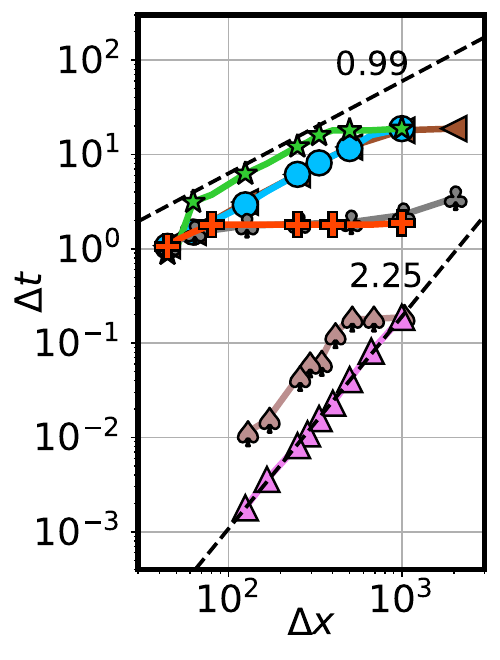} &
     \raisebox{-0.16cm}{\includegraphics[width=0.271\linewidth]{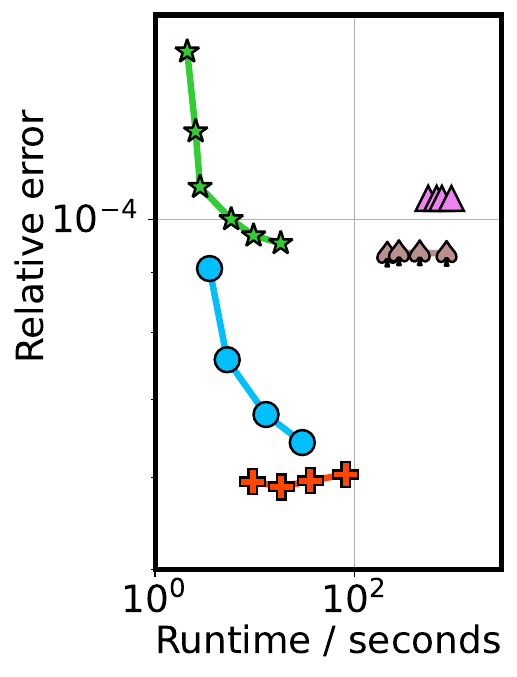}} &    
    \hspace{-0.1cm}\raisebox{2.5cm}{\includegraphics[width=0.21\linewidth]{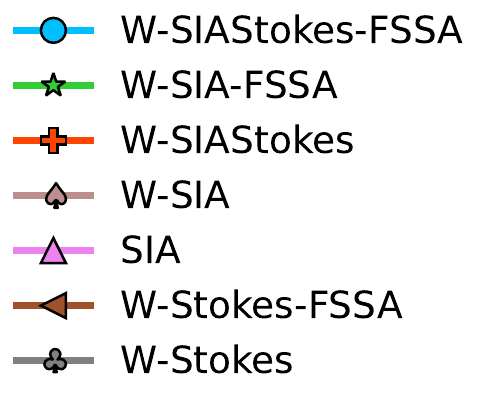}}
\end{tabular}
\caption{Left plot: scaling of the largest feasible time step $\Delta t_*$ as a function of the horizontal mesh size $\Delta x$ when the vertical mesh size 
$\Delta y = 90$ meters is fixed.  
Right plot: the model error as a function of the wall clock runtime when the nonlinear Stokes solution is taken as a reference. 
In all cases the final time is fixed at $t=20$ years. Mesh sizes $\Delta x = 250$ meters and $\Delta y = 90$ meters are also fixed. 
The time step is refined starting at the formulation largest feasible time step $\Delta t = \Delta t_*, \Delta t_*/2, \Delta t_*/4,...$
The time step for the reference solution is fixed at $\Delta t = 0.1$ years.}
\label{fig:slabslope:dt_vs_dx_and_runtime}
\end{figure}

\subsubsection{Time step restriction scaling}
Now we compute $\Delta t_*$ as a function of $\Delta x$ as described in Section \ref{sec:computations_largest_timestep}.
We run the simulations with and without the FSSA terms defined in the scope of Section \ref{sec:fssa}. 
The final simulation times are adjusted to the magnitude of the largest time step sizes making the comparison more realistic. We use the final simulation times: 
$t=100$ years (W-SIAStokes-FSSA and W-SIA-FSSA), $t=12$ years (W-SIAStokes and W-SIA), $t=5$ years (SIA).
The results are shown in the first plot of Figure \ref{fig:slabslope:dt_vs_dx_and_runtime}. We observe that the largest stable time step 
$\Delta t_*$ is allowed to be from 10 times (coarse resolution) to 100 times (fine resolution) larger when using W-SIAStokes as compared to W-SIA.
Furthermore, in W-SIAStokes, $\Delta t_*$ has a constant relation to $\Delta x$ over the whole range of the chosen $\Delta x$ except for the finest resolution 
where the relation becomes linear. The allowed time steps when using SIA and W-SIAStokes are small. 
The $\Delta_* t$ vs. $\Delta x$ scaling in the latter two formulations is quadratic which is less desired.

\subsubsection{Run-time versus accuracy}
Next, we perform an experiment measuring the ratio between the relative model error and the computational (wall clock) runtime for each of the SIA formulations. 
This is important as a small $\Delta t$ does not necessarily imply a small computational time for the whole simulation. This depends on the computational time required to evaluate 
the velocity functions in each time-step. 
We set the final simulation 
time to $t=20$ years. 
For W-SIAStokes-FSSA we used $\Delta t = 6, 3, 1.5, 0.75$. 
For W-SIA-FSSA we used $\Delta t = 12, 6, 3, 1.5, 0.75, 0.4$. 
For W-SIAStokes we used $\Delta t = 1.8, 0.9, 0.45, 0.2$. 
For W-SIA we used $\Delta t = 0.04, 0.03, 0.02, 0.01$, and for SIA we used $\Delta t = 0.008, 0.007, 0.006, 0.005$. 
The error is computed with the nonlinear Stokes solution as a reference with $\Delta t = 0.1$ years. 
In all cases the mesh sizes $\Delta x = 250$ meters and $\Delta y = 90$ meters 
are fixed. 
The result is given in the 
second plot of Figure \ref{fig:slabslope:dt_vs_dx_and_runtime}. We observe that W-SIAStokes-FSSA allows small computational runtimes with only a mild increase of the model error. 
W-SIA-FSSA also allows for small computational runtimes, but the model error is larger than in the W-SIAStokes-FSSA case. 
The unstabilized (weak nor strong) SIA formulations do not allow for small computational runtimes, except in the case when we use the 
W-SIAStokes formulation.
\begin{figure}
    \centering
\begin{tabular}{ccc}
    \multicolumn{3}{c}{\textbf{Slab on a slope surface case, varying $\mathbf{\theta}$ (FSSA parameter), W-SIAStokes-FSSA}} \vspace{0.2cm} \\
    \hspace{1cm}$\mathbf{\Delta y = 83.3} \textbf{ meters}$ & \hspace{1cm}$\mathbf{\Delta y = 41.67} \textbf{ meters}$ & \\    
    \includegraphics[width=0.28\linewidth]{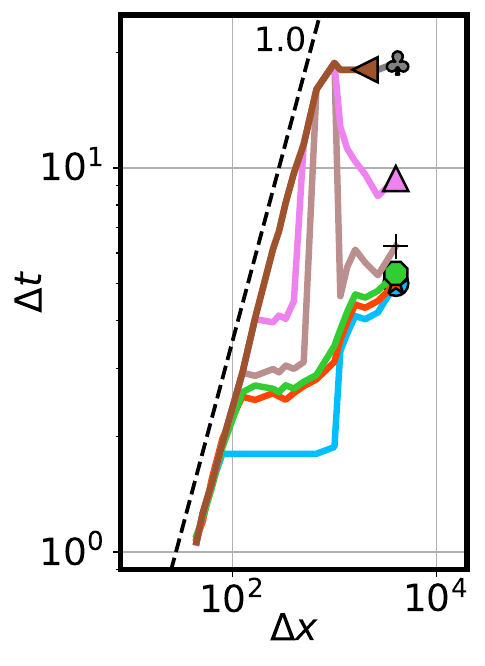} &
    \raisebox{-0.01cm}{\includegraphics[width=0.28\linewidth]{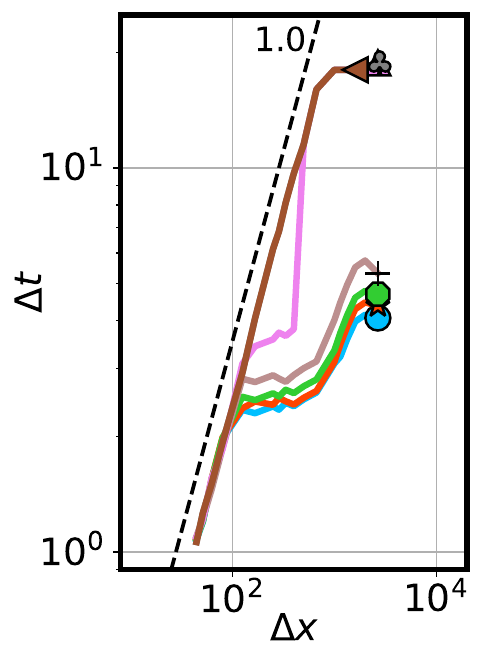}} &
    \raisebox{1.2cm}{\hspace{0.3cm}\includegraphics[width=0.17\linewidth]{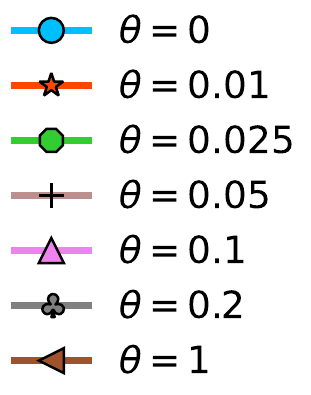}}
\end{tabular}
\caption{Largest feasible time step $\Delta t_*$ as a function of the horizontal mesh size $\Delta x$ for the W-SIAStokes-FSSA formulation, when the 
FSSA parameter $\theta$ varies. Vertical mesh size is fixed at
$\Delta y = 83.3$ meters (left plot), and at $\Delta y = 41.67$ meters.}
\label{fig:slabslope:dt_vs_dx_theta_fssa}
\end{figure}

\subsubsection{Impact of the FSSA parameter $\theta$ on the time step restriction scaling}
As we demonstrated in the previous paragraphs W-SIAStokes-FSSA allows for largest time steps among 
all the considered formulations. The FSSA parameter in those experiments was fixed at $\theta = 1$. 
In Figure \ref{fig:slabslope:dt_vs_dx_theta_fssa} we illustrate the effect of the choice of $\theta$ 
on the scaling of the largest feasible time step $\Delta t_*$ as a function of $\Delta x$. We observe that choice of small $\theta$ leads 
to a non-robust $\Delta t_*$ vs. $\Delta x$ scaling behavior, where $\Delta t_*$ has to be taken around $10$-times smaller as $\theta \to 0$. 
As the FSSA parameter approaches $\theta = 0.2$ then the $\Delta t_*$ vs. $\Delta x$ 
scaling approaches the linear behavior as also observed in Figure \ref{fig:slabslope:dt_vs_dx_and_runtime}. Furthermore the scaling stays the same for $\theta > 0.2$. 
This implies that, for the given test case, the choice of $\theta$ is not sensitive, and can be left at $\theta = 1$ without losing the length of the largest 
feasible time step.


\subsection{An idealized ice sheet surface case}
\label{sec:experiments:vialov}
\subsubsection{Configuration}
We now change configuration in order to study the impact of ice sheet margins. The ice sheet configuration in this test case takes horizontal values $x \in [-L, L]$ and vertical values $y \in [0, H]$, 
where the ice sheet half length 
is $L = 750 \cdot 10^{3}$ meters and the ice sheet height is $H = 3 \cdot 10^{3}$ meters. 
To construct the ice sheet surface for this test case we define the following auxilary profile:
$$h_1(x) = \left(3 - \left(\frac{x}{L}\right)^2 \right)^{0.58},$$
and then use it in the initial surface definition:
$$h (x,0) = H\, \frac{h_1(x) - h_1(-L)}{h_1(0)}.$$

When solving the free-surface equation \eqref{eq:free_surface_eq} we set Dirichlet boundary conditions $h(-L,t) = h(-L,0)$ on the left 
boundary, and $h(L,t) = h(L,0)$ on the right boundary. 
When solving one of the momentum balances we impose no-slip boundary conditions on $\Gamma_b$ and on $\Gamma_l$, whereas on $\Gamma_s$ 
we set stress-free boundary conditions.

The relation between the number of horizontal mesh elements $n_x$ and horizontal mesh size $\Delta x$ 
that we use to perform the experiments is given in the table below. 

\vspace{0.2cm}
\begin{center}
\begin{tabular}{c|cccccc}
    \multicolumn{7}{c}{\textbf{An idealized ice sheet surface case}} \vspace{0.1cm} \\
    \multicolumn{7}{c}{\textbf{Number of elements and the mesh size}} \vspace{0.2cm} \\
    $\mathbf{n_x}$ & 80 & 200 & 400 & 600 & 800 & 1000 \\ \hline
$\mathbf{\Delta x}$ & 18750 & 7500 & 3750 & 2500 & 1875 & 1500 \vspace{0.3cm}
\end{tabular}
\end{center}
\vspace{0.3cm}

Note that the number of vertical mesh elements is $n_y=12$ (corresponding to vertical mesh size $\Delta y = 250$ meters) is fixed for all experiments unless stated otherwise. The FSSA scaling parameter is always fixed at $\theta = 1$.

\subsubsection{Time step restriction scaling}
In the left plot of Figure \ref{fig:vialov:dt_vs_dx_and_runtime} we display the largest feasible time step $\Delta t_*$ as a function 
of $\Delta x$ for the different SIA formulations. Time step $\Delta t_*$ scales linearly with respect to $\Delta x$ for W-SIAStokes 
and W-SIAStokes-FSSA. 
For SIA and W-SIA 
we observe a quadratic scaling. The benefits when using the FSSA stabilization terms together with W-SIA 
disappear as the $\Delta x$ is fine enough. Among all the formulations the largest time steps can be taken when W-SIAStokes W-SIAStokes-FSSA are used. 
For W-SIAStokes-FSSA the time steps increase 
approximately 100 times, across the whole range of the considered mesh sizes.

\subsubsection{Run-time versus accuracy}
\begin{figure}[h!]
    \centering
    \begin{tabular}{ccc}
    \multicolumn{3}{c}{\textbf{\hspace{-2cm}An idealized ice sheet surface case}} \vspace{0.2cm} \\
    \hspace{0.8cm}\textbf{$\mathbf{\Delta t}$ vs. $\mathbf{\Delta x}$} & \hspace{0.8cm}\textbf{Error vs. Runtime} & \\    
     \includegraphics[width=0.25\linewidth]{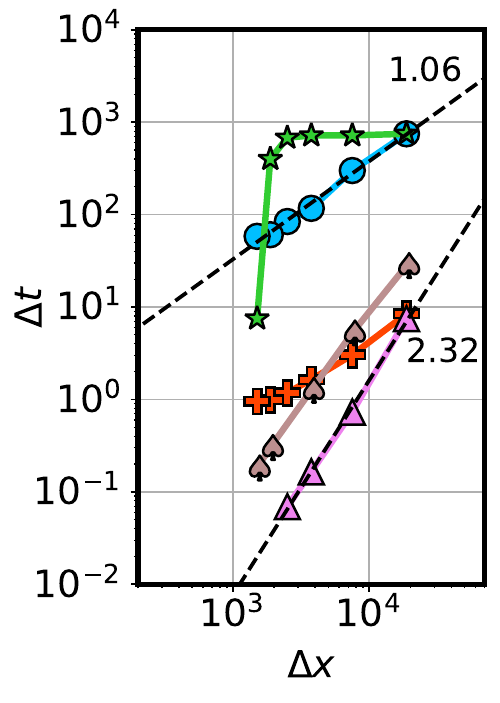} &
     \raisebox{0.06cm}{\includegraphics[width=0.244\linewidth]{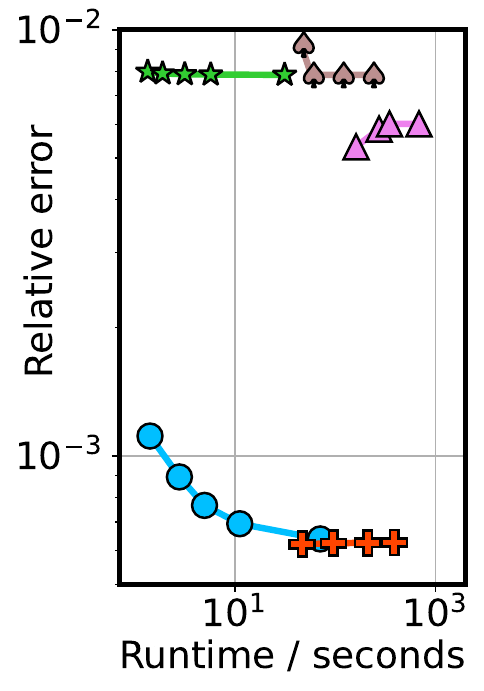}} &    
    \hspace{-0.3cm}\raisebox{3.5cm}{\includegraphics[width=0.17\linewidth]{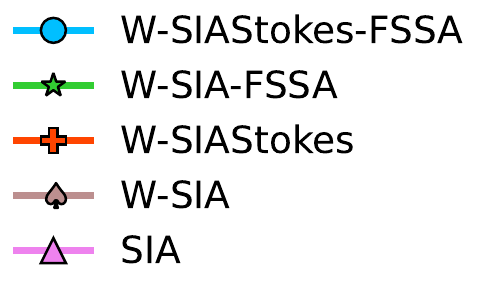}}
\end{tabular}
\caption{Left plot: scaling of the largest feasible time step $\Delta t_*$ as a function of the horizontal mesh size $\Delta x$ when the vertical mesh size 
$\Delta y = 90$ meters is fixed. 
Right plot: the model error as a function of the wall clock runtime when the nonlinear Stokes solution is taken as a reference. 
In all cases the final time is fixed at $t=20$ years. Mesh sizes $\Delta x = 250$ meters and $\Delta y = 90$ meters are also fixed.}
\label{fig:vialov:dt_vs_dx_and_runtime}
\end{figure}

In the next experiment we compute the ratio between the computatinal runtime (wall clock) and the modeling error. 
For all formulations we fixed the horizontal mesh size to $\Delta x = 3750$ meters, and ran the simulation until $t=100$ years. 
The error is computed using the nonlinear Stokes solution as a reference, where the time step is $\Delta t = 0.1$ years. 
The maximum time steps for testing each SIA formulation is taken in line 
with the largest feasible time step for $\Delta x = 3750$ meters from Figure \ref{fig:vialov:dt_vs_dx_and_runtime}. 
For both FSSA stabilized weak SIA formulations we used $\Delta t = 50, 25, 12.5, 6, 1$. 
For W-SIAStokes we used $\Delta t = 1.6, 0.8, 0.4, 0.2$. 
For W-SIA we used $\Delta t = 1, 0.8, 0.4, 0.2$, and for SIA 
we used $\Delta t = 0.16, 0.1, 0.08, 0.04$. 
The results are presented in Figure \ref{fig:vialov:dt_vs_dx_and_runtime}, plot on the right.
Runtime vs. error ratio of W-SIASTokes is favorable over all the other 
formulations that we consider. This is, both, when the FSSA terms are present, and when the FSSA terms are not present.

\subsubsection{Long time simulations}
\begin{figure}[h!]
    \centering
    \begin{tabular}{cc}
    \multicolumn{2}{c}{\textbf{\hspace{-0cm}Idealized ice sheet surface case}} \vspace{0.2cm} \\
    \multicolumn{2}{c}{\hspace{0cm}\textbf{Surface elevations at $\mathbf{t=10^4}$ years}} \vspace{0.2cm}\\
    \textbf{W-SIAStokes-FSSA} & \textbf{W-SIA-FSSA} \\
    \includegraphics[width=0.5\linewidth]{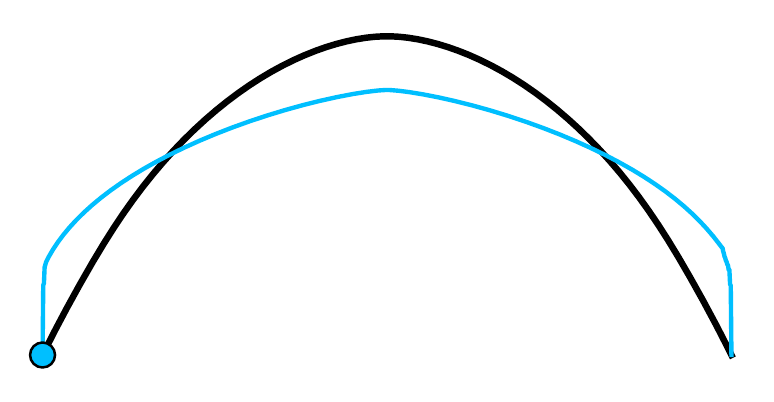} & 
    \includegraphics[width=0.5\linewidth]{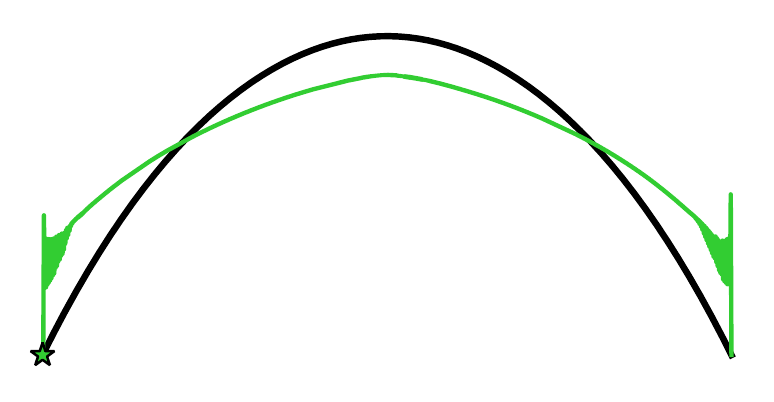}
\end{tabular}
\caption{Two surface evolutions after $t=10^4$ years. The horizontal mesh size is $\Delta x = 1500$ meters, whereas the vertical mesh size is $\Delta y = 250$ meters. 
The time steps take the largest feasible values for the given mesh sizes: $\Delta t = 58.2$ years for W-SIAStokes, 
and $\Delta t = 7.5$ years for W-SIA. Both formulations are stabilized by means of the FSSA stabilization terms with the FSSA parameter set to $\theta = 1$.}
\label{fig:vialov:surface_solutions}
\end{figure}

We now compute the surface evolution over a long period of time: the final time is $t=10^4$ years. 
We choose a fine mesh size $\Delta x = 1500$ meters. The time step is chosen 
in line with Figure \ref{fig:vialov:dt_vs_dx_and_runtime} (left plot) for the given $\Delta x$, that is, $\Delta t = 58.2$ years for W-SIAStokes, 
and $\Delta t = 7.5$ years for W-SIA. 
In Figure \ref{fig:vialov:surface_solutions} we show the solutions obtained from the two formulations. 
The solution obtained using W-SIAStokes does not entail any oscillations, whereas the solution in the W-SIA case 
entails oscillations close to the lateral boundaries. 
We have not fully explored the behavior. However, we speculate that this is due to the lack of the control over all strain components 
in W-SIA (see \ref{eq:sia_weakform_less}) in contrast with 
W-SIAStokes \eqref{eq:sia_weakform_all}, where all the strain components are present. 
This implies that a bound $\|\bm D \bm u\|_{L^2(\Omega)} \leq \| \bm f\|_{L^2(\Omega)}$, where $\bm f$ is the gravity field, 
can be obtained by using Korn's inequality. 
This provides control over the derivatives of the velocities $\bm u$, which prevents $\bm u$ from being too oscillatory. 

\subsection{A two-dimensional cross-section of Greenland}
\label{sec:experiments:greenland}
\subsubsection{Configuration}
\begin{figure}[h!]
    \centering
\begin{tabular}{cc}
\textbf{Greenland} & \textbf{The 2D Greenland geometry}\vspace{0.2cm} \\
    \includegraphics[width=0.2\linewidth]{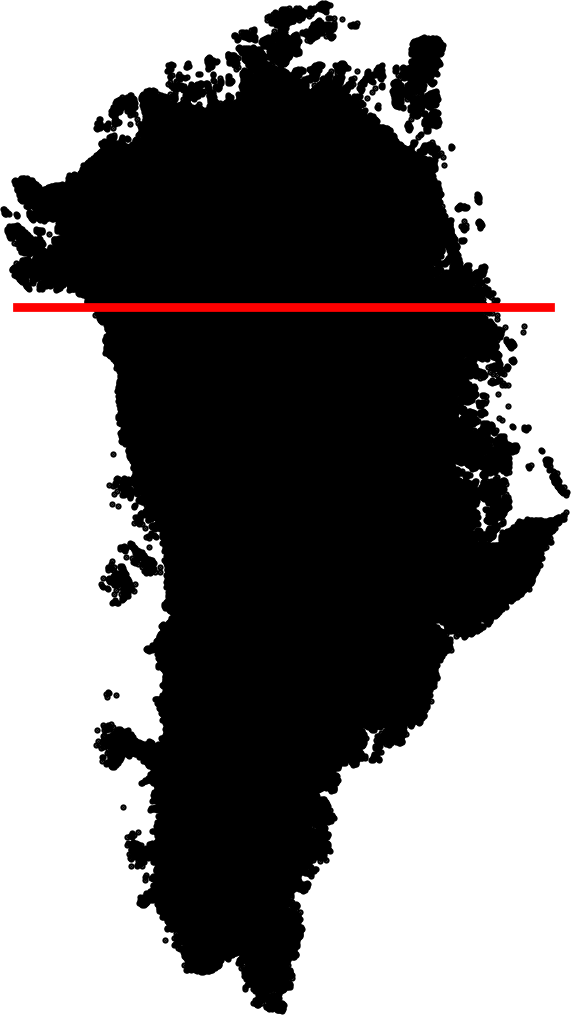} &
\raisebox{1.45cm}{\includegraphics[width=0.4\linewidth]{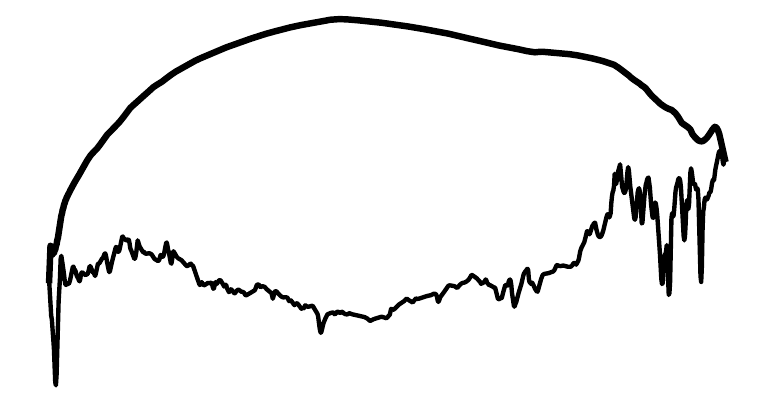}}
\end{tabular}
\caption{Top view over a three-dimensional Greenland geometry \cite{bedmachine} (left), where the intersecting 
red line gives the boundary of the Greenland cross-section (right) used as one of the computational domains in this paper.}
\label{fig:greenland_line_crossection}
\end{figure}

In this test case we still consider a two-dimensional ice sheet geometry, but 
with a more realistic initial surface elevation as well as a more realistic bedrock elevation. 
We simplify the full three-dimensional Greenland geometry obtained from BedMachine \cite{bedmachine}, 
and intersect it with a horizontal line to get the boundary points over a cross section as displayed in Figure 
\ref{fig:greenland_line_crossection}. Then we represent the surface elevation and the bedrock elevation by using a cubic spline interpolation. 
This allows for evaluating the surfaces at an arbitrary location, which we employ when considering horizontally varying 
mesh sizes in our computational study.

When solving the free-surface equation \eqref{eq:free_surface_eq} we set Dirichlet boundary conditions $h(-L,t) = h(-L,0)$ on the left 
boundary, and $h(L,t) = h(L,0)$ on the right boundary. 
When solving one of the momentum balances we impose no-slip boundary conditions on $\Gamma_b$ and on $\Gamma_l$, whereas on $\Gamma_s$ 
we set stress-free boundary conditions.

The relation between the number of horizontal mesh elements $n_x$ and horizontal mesh size $\Delta x$ 
that we use to perform the experiments is given in the table below. 

\vspace{0.2cm}
\begin{center}
\begin{tabular}{c|cccccc}
    \multicolumn{7}{c}{\textbf{Greenland (2D) profile case}} \vspace{0.1cm} \\
    \multicolumn{7}{c}{\textbf{Number of elements and the mesh size}} \vspace{0.2cm} \\
    $\mathbf{n_x}$ & 80 & 200 & 400 & 600 & 800 & 1000 \\ \hline
$\mathbf{\Delta x}$ & 12362 & 4945 & 2472 & 1648 & 1236 & 989 \vspace{0.3cm}
\end{tabular}
\end{center}
\vspace{0.3cm}

\subsubsection{Time step restriction scaling}
\begin{figure}
    \centering
\begin{tabular}{ccc}
    \multicolumn{3}{c}{\hspace{-2cm}\textbf{Greenland (2D) surface case}} \vspace{0.2cm}\\
    \hspace{0.8cm}\textbf{$\mathbf{\Delta t}$ vs. $\mathbf{\Delta x}$} & \hspace{0.8cm}\textbf{$\mathbf{\Delta t}$ vs. $\mathbf{\Delta x}$} & \\    
    \hspace{0.8cm}(no upwind viscosity) & \hspace{0.65cm}(with upwind viscosity) & \\
    \includegraphics[width=0.25\linewidth]{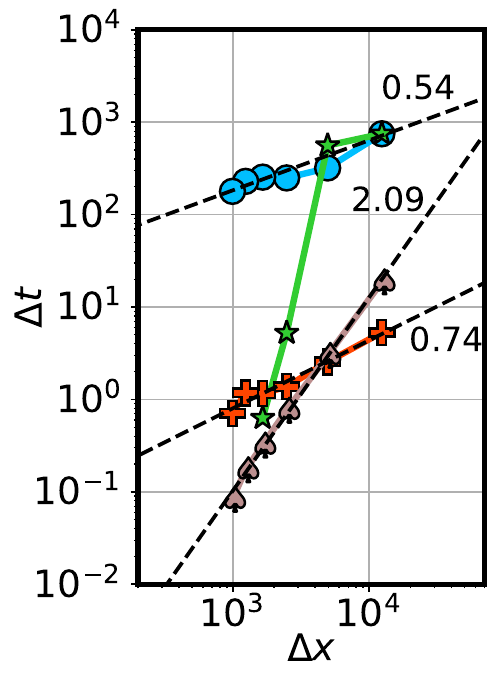} &
     \includegraphics[width=0.25\linewidth]{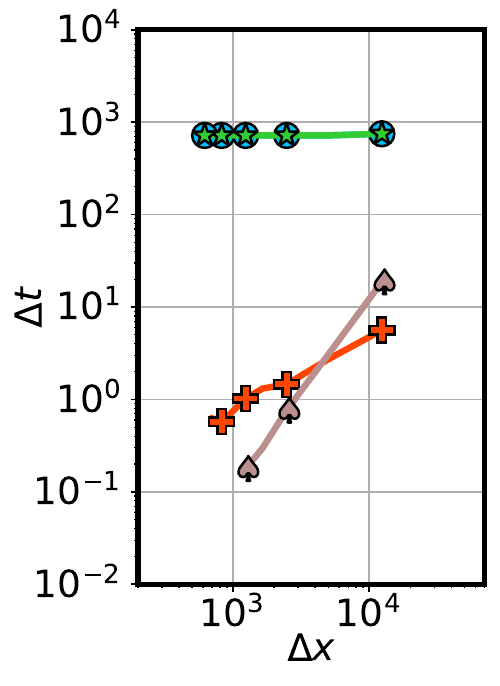} &
    \hspace{-0.5cm}\raisebox{3.5cm}{\includegraphics[width=0.17\linewidth]{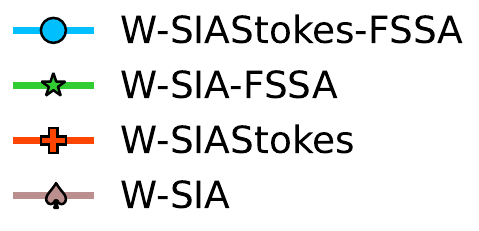}}
\end{tabular}
\caption{
    Left plot: scaling of the largest feasible time step $\Delta t_*$ as a function of the horizontal mesh size $\Delta x$ when the vertical mesh size 
is fixed. 
Right plot: the model error as a function of the wall clock runtime when the nonlinear Stokes solution is taken as a reference. 
In all cases the final time is fixed at $t=400$ years. The horizontal mesh size is fixed at $\Delta x=2472$ meters.}
\label{fig:greenland:dt_vs_dx}
\end{figure}

In Figure \ref{fig:greenland:dt_vs_dx} (left plot) we investigate 
the scaling of $\Delta t$ as a function of $\Delta x$. 
We observe that in the case of W-SIAStokes-FSSA the order of scaling is almost $1$, 
whereas in the W-SIAStokes-FSSA case the order of scaling is close to $0.5$.
The scaling in the W-SIA case is approximately of order $2$, which is similar as in all previous test cases. 
When W-SIA-FSSA is used the scaling behaves unpredictably, 
similar as in Section \ref{sec:experiments:vialov}.
Among all the formulations, 
the time steps are the largest in the W-SIAStokes-FSSA.

\subsubsection{Upwinding}
After visualizing the computed solutions using all the considered SIA formulations we observed 
spurious oscillations over the western part of the two-dimensional Greenland geometry (see Figure \ref{fig:greenland:solutions_surface}). 
For that reason 
we in addition combined the SIA formulations with the first-order viscosity (also upwind viscosity) operator 
added to the free-surface equation. The role of that operator is dampening of 
the spurious oscillations. In Figure \ref{fig:greenland:dt_vs_dx} (right plot) we compute 
the scaling of $\Delta t$ as a function of $\Delta x$ when the first-order viscosity 
operator is added to the free-surface equation, for each of the considered SIA formulations. We 
observe that the scaling does not change when no FSSA terms are added W-SIAStokes or W-SIA. 
However, the scaling, when the FSSA terms are employed, changes from linear (no first-order viscosity) 
to constant (with first-order viscosity). The latter is favorable. 

\begin{figure}
    \centering
\begin{tabular}{cc}
    \multicolumn{2}{c}{\hspace{0cm}\textbf{Greenland (2D) surface case}} \vspace{0.2cm}\\
    \multicolumn{2}{c}{\hspace{0cm}\textbf{Surface elevation comparison against the reference at $\mathbf{t=400}$ years}} \vspace{0.2cm}\\
    \textbf{W-SIAStokes-FSSA} & \textbf{W-SIA-FSSA}\\
     \includegraphics[width=0.5\linewidth]{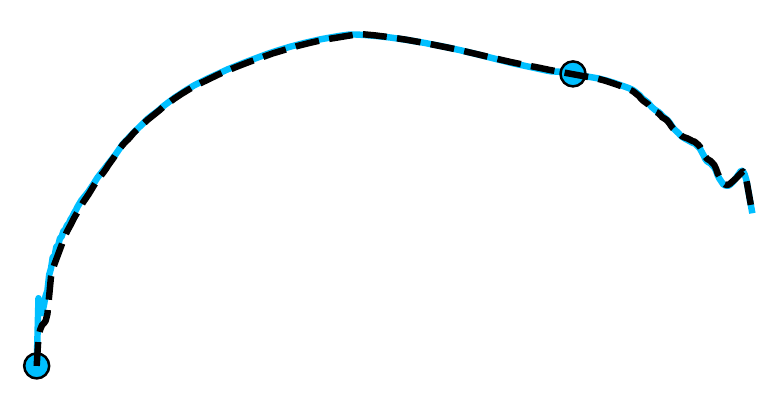} &
     \includegraphics[width=0.5\linewidth]{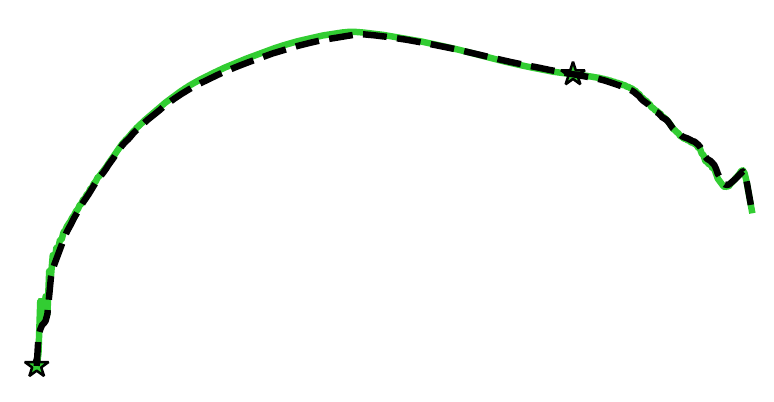} \\
     \textbf{W-SIAStokes} & \textbf{W-SIA}\\
     \includegraphics[width=0.5\linewidth]{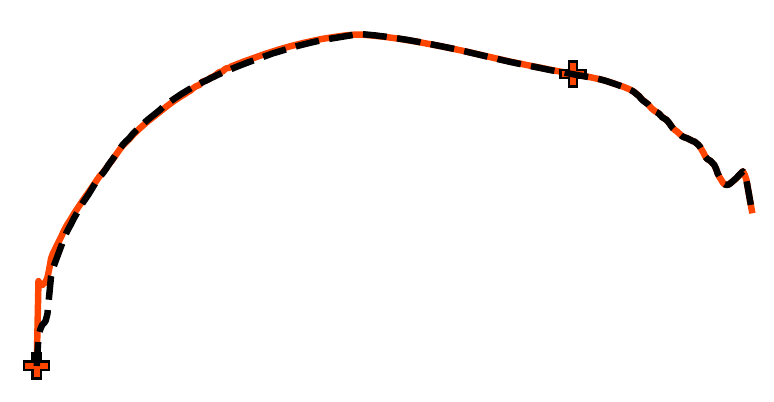} &
     \includegraphics[width=0.5\linewidth]{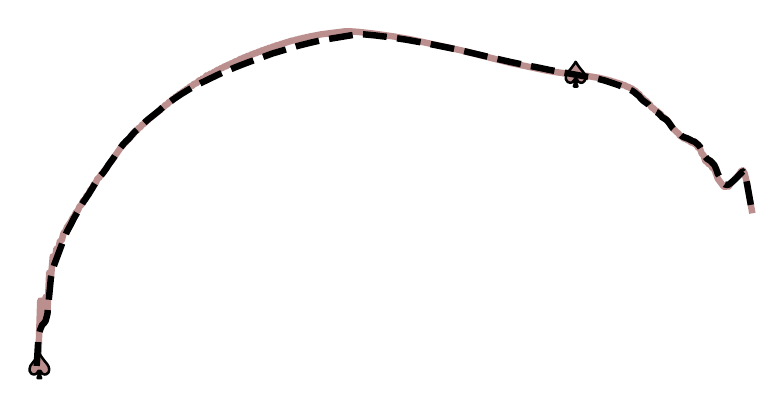} \\
     \textbf{W-SIAStokes-FSSA-UV} & \textbf{W-SIA-FSSA-UV}\\
     \includegraphics[width=0.5\linewidth]{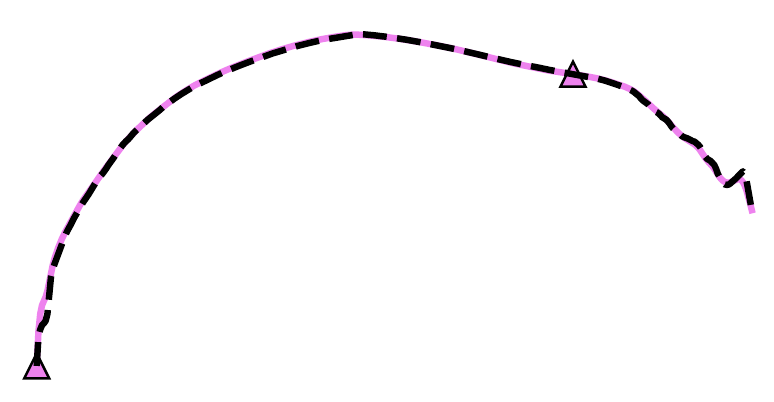} &
     \includegraphics[width=0.5\linewidth]{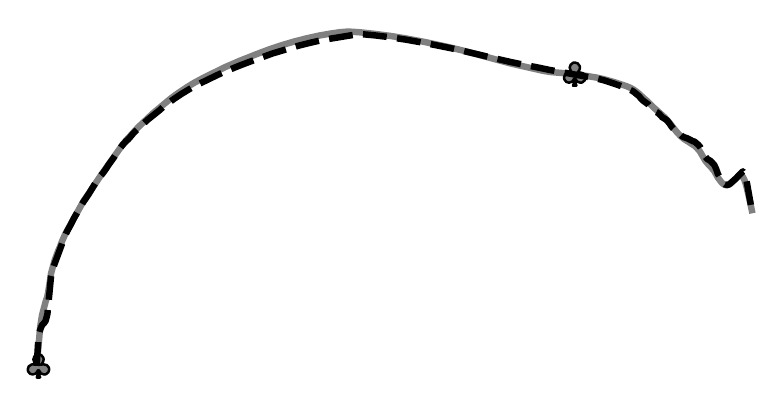}
\end{tabular}
\caption{Surface evolutions computed using different SIA formulations, 
after $t=400$ years. The horizontal mesh size is $\Delta x = 2472$ meters.
The time steps take the largest feasible values for the given mesh sizes (see Figure \ref{fig:greenland:dt_vs_dx}).}
\label{fig:greenland:solutions_surface}
\end{figure}

\subsubsection{Long time simulations}
In Figure \ref{fig:greenland:solutions_surface} we, for each of the considered SIA formulations, 
plot the surface elevations after $t=400$ years of simulation time, 
at $\Delta x = 2472$ meters, and make a comparison towards the surface 
elevations computed using the nonlinear Stokes problem (reference). 
We observe that all the SIA formulation solutions are overall 
a good approximation to 
the reference solution. As stated in the previous paragraph, all solutions in the first three rows of Figure \ref{fig:greenland:solutions_surface} 
involve spurious oscillation 
at the western part of the two-dimensional Greenland geometry.
The plots in the last row of Figure \ref{fig:greenland:solutions_surface} 
display the viscous solutions, and we observe that the oscillations have disappeared.

\begin{figure}
    \centering
\begin{tabular}{cc}
    \multicolumn{2}{c}{\hspace{-2cm}\textbf{Greenland (2D) surface case}} \vspace{0.2cm}\\
    \hspace{0.8cm}\textbf{Error vs. Runtime}  & \\    
     \includegraphics[width=0.3\linewidth]{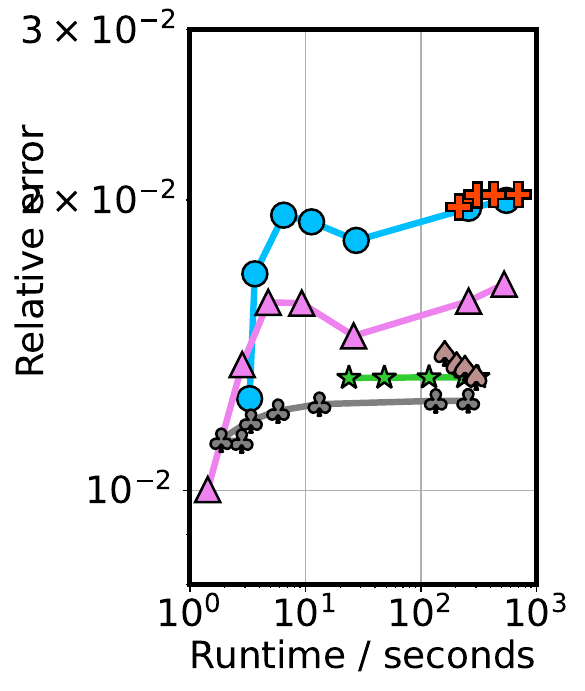} &
    \hspace{-0.5cm}\raisebox{2.5cm}{\includegraphics[width=0.2\linewidth]{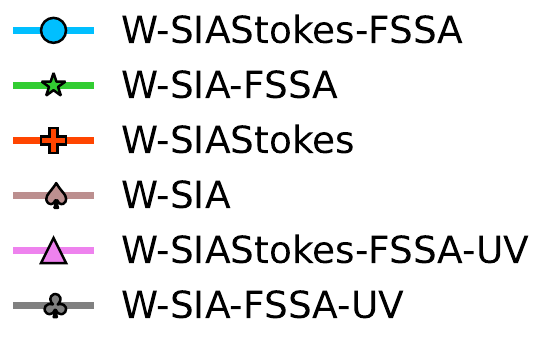}}
\end{tabular}
\caption{Approximate model error as a function of the wall clock runtime when the nonlinear Stokes solution is taken as a reference. 
In all cases the final time is fixed at $t=400$ years. The horizontal mesh size is $\Delta x = 2472$.}
\label{fig:greenland:error_vs_runtime}
\end{figure}

\subsubsection{Run-time versus accuracy}
In Figure \ref{fig:greenland:error_vs_runtime} 
compute the ratio between the computatinal runtime (wall clock) and the modeling error. 
For all formulations we fixed the horizontal mesh size to $\Delta x = 2472$ 
meters, and ran the simulation until $t=400$ years. 
The error is computed using the nonlinear Stokes solution as a reference, where the time step is $\Delta t = 0.4$ years. 
The maximum time steps for testing each SIA formulation is taken in line 
with the largest feasible time step for $\Delta x = 2472$ meters from 
Figure \ref{fig:greenland:dt_vs_dx}. 
For W-SIAStokes-FSSA we used $\Delta t = 249, 125, 60, 30, 10, 1, 0.5$ years. 
For W-SIA-FSSA we used $\Delta t = 5, 2.5, 1, 0.5, 0.4$ years. 
For W-SIAStokes we used $\Delta t = 1.37, 0.9, 0.65, 0.4$ years. 
For W-SIA we used $\Delta t = 0.76, 0.6, 0.5, 0.4$. 
When W-SIAStokes-FSSA and W-SIA-FSSA are further augmented with 
the with upwind viscosity in the free-surface equation 
we used $\Delta t = 800, 400, 200, 100, 50, 25, 12.5$ years. 
We observe that the FSSA stabilized weak formulations augmented with the upwind viscosity term have by far the best 
error vs. runtime ratio.

\begin{figure}
    \centering
\begin{tabular}{cc}
    \multicolumn{2}{c}{\hspace{0cm}\textbf{Greenland (2D) surface case}} \vspace{0.2cm}\\
    \multicolumn{2}{c}{\hspace{0cm}\textbf{Surface elevations at $\mathbf{t=10000}$ years}} \vspace{0.2cm}\\
    \textbf{W-SIAStokes-FSSA} & \textbf{W-SIA-FSSA}\\
     \includegraphics[width=0.5\linewidth]{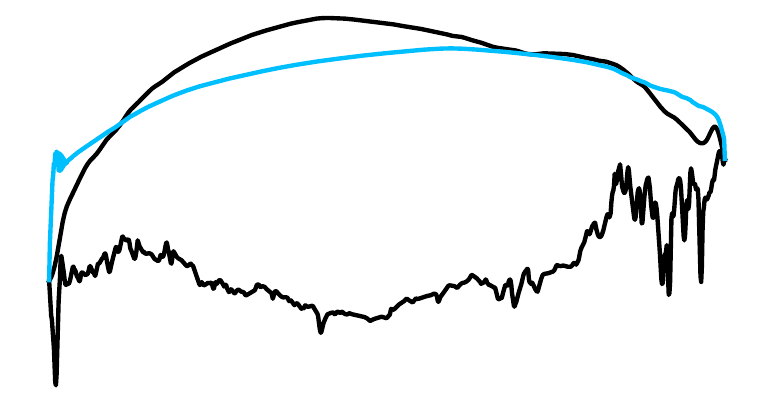} &
     \includegraphics[width=0.5\linewidth]{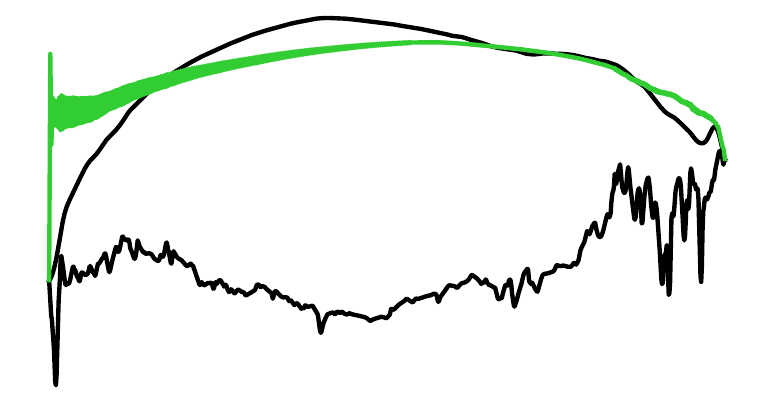} \\
     \textbf{W-SIAStokes-FSSA-UV} & \textbf{W-SIA-FSSA-UV}\\
     \includegraphics[width=0.5\linewidth]{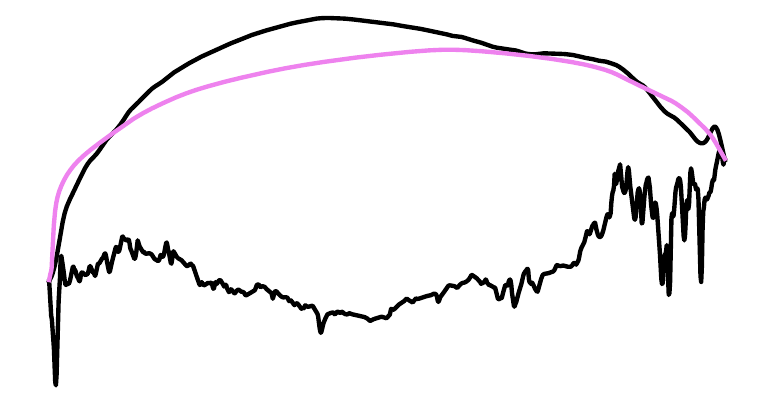} &
     \includegraphics[width=0.5\linewidth]{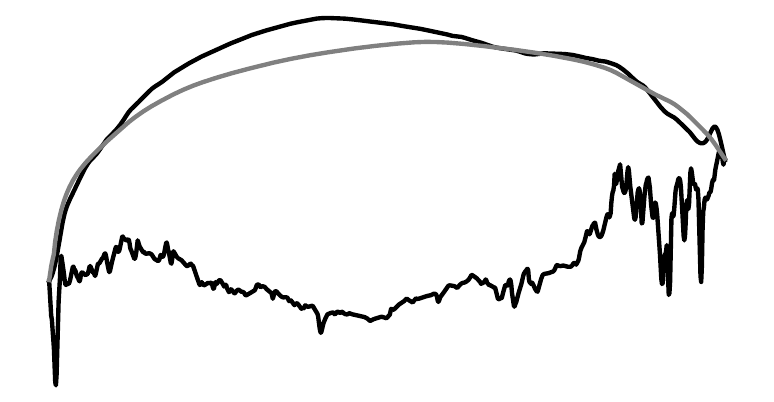}
\end{tabular}
\caption{Surface elevations computed using the FSSA stabilized SIA formulations, with and without 
the first-order (upwind) viscosity added to the free-surface equation. 
The surface elevations are evaluated at $t=1000$ years. 
The horizontal mesh size is $\Delta x = 2472$ meters.
The time steps take the largest feasible values for the given mesh sizes (see Figure \ref{fig:greenland:dt_vs_dx}).}
\label{fig:greenland:solutions_surface}
\end{figure}
\subsection{Computational cost estimates with the parameters inferred from the numerical experiments}
\label{section:performance_analysis_revisited}
In the subsections above we gathered experimental data on the slope of the scaling of the time step restriction as a function 
of the horizontal mesh size. In this subsection we use those slopes in place of the parameter $\gamma$ 
in \ref{section:performance_analysis}, 
to speculate on the relation between computational work and the number of the horizontal mesh nodes  $m$ 
for the different momentum models. 
Across the experiments we observe that W-SIA-FSSA, W-SIAStokes, and W-SIAStokes-FSSA have a linear time-step constraint $\gamma=1$. 
On Greenland with upwind viscosity we saw that $\gamma$ can also be smaller than $1$, 
but we here chose the worst case scenario value ($\gamma=1$) for W-SIA-FSSA and W-SIAStokes. 
In the W-SIA model case the worst case is $\gamma=2$. The same holds for the standard SIA model case. 
We gather these results in the table below. In the W-Stokes-FSSA case we get $\gamma=1$, wherease in the W-Stokes case 
we have $\gamma=1$ when the horizontal mesh size is small, and $\gamma=0$ for larger horizontal mesh sizes -- the worst case is then $\gamma=1$.
%
\\
\begin{center}
{\tabulinesep=0.6mm
\begin{tabu}{ l || c | c | c  }
  \textbf{Model} & $\bm \gamma$ &  \textbf{Computational cost estimate} & \textbf{Evaluated with } $\bm\gamma$ \textbf{and} $\mathbf{d=3}$\\ \hline \hline 
  \textbf{W-Stokes} & $1$ & $C_S\, m^{1+\gamma/(d-1)+\alpha}$ &  $C_S\, m^{1.5+\alpha}$ \\ \hline
  \textbf{W-Stokes-FSSA} & $1$ & $C_S\, m^{1+\gamma/(d-1)+\alpha}$ & $C_S\, m^{1.5+\alpha}$  \\ \hline
  \textbf{W-SIAStokes} & $1$ & $\frac{1}{N_{\text{iter}}}\, C_{S}\, m^{1+\gamma/(d-1)+\alpha}$ &  $\frac{1}{N_{\text{iter}}}\, C_{S}\, m^{1.5+\alpha}$ \\ \hline
  \textbf{W-SIAStokes-FSSA} & $1$ & $\frac{1}{N_{\text{iter}}}\, C_{S}\, m^{1+\gamma/(d-1)+\alpha}$ & $\frac{1}{N_{\text{iter}}}\, C_{S}\, m^{1.5+\alpha}$ \\ \hline
  \textbf{W-SIA} & $2$ & $\frac{d+1}{(d+1)^{1+\alpha}}\, \frac{1}{N_{\text{iter}}}\, C_{S}\, m^{1+\gamma/(d-1)+\alpha}$ & $\frac{4}{4^{1+\alpha}}\, \frac{1}{N_{\text{iter}}}\, C_{S}\, m^{2+\alpha}$  \\ \hline
  \textbf{W-SIA-FSSA} & $1$ & $\frac{1}{N_{\text{iter}}}\, C_{S}\, m^{1+\gamma/(d-1)+\alpha}$ &  $\frac{1}{N_{\text{iter}}}\, C_{S}\, m^{1.5+\alpha}$\\ \hline
  \textbf{SIA} & $2$ & $C_{\text{SIA}}\, m^{1+\gamma/(d-1)}$ &  $C_{\text{SIA}}\, m^{2}$
\end{tabu}
}
\end{center}
We observe that the SIA model has the lowest asmyptotic cost scaling $m^2$. The next best is W-SIAStokes-FSSA with $m^{2.5}$ when $\alpha=1$ (sparse direct solver). 
We argue that the computational cost is in the two cases comparable when it comes to the asmyptotic scaling as $m\to\infty$. However, 
knowing that the SIA model is a crude simplification of the Stokes model \eqref{eq:fullstokes}, and that the W-SIAStokes model 
only carries minor simplified elements (the viscosity) as compared to W-Stokes, 
the W-SIAStokes model is more accurate as we could also observe from the runtime vs. error experiments figures. 
The computational cost in the W-SIA model also scales as $m^{1.5+\alpha}$, but this does not hold true for all the horizontal mesh size choices 
as observed from the experiments.
We note that $N_{nlin}$ is large for large $\Delta t$. 
This is due to that the nonlinear iteration initial guess to compute the solution at time $t_{k+1}$ is the solution from time $t_k$, and the two 
solutions when $\Delta t = t_{k+1} - t_k$ is large, are typically very different, implying that the nonlinear iteration count is large.

The time step sizes can also be restricted by some other components 
of an ice sheet model, 
such as the temperature evolution or climate data. 
However, when these time step sizes are small, it is still not certain that the velocity solution has to be updated at the same small time 
time-step as e.g. the temperature. This study is out of the scope of this paper.

\section{Final remarks}
\label{section:conclusion}
In this paper we investigated the benefits when using the SIA momentum balance model 
\eqref{eq:SIAeqa} written on different weak formulations. We referred to the weak SIA model \eqref{eq:sia_weakform_less} as W-SIA, 
the extended weak SIA model alias weak form linear Stokes equations employing the SIA viscosity function as W-SIAStokes, 
the weak full nonlinear Stokes problem as W-SIA, and to the standard (strong form) SIA model as SIA. When the different SIA forms 
were additionally stabilized by mean of the FSSA terms \eqref{eq:fssa} we appended an abbreviation FSSA to each of the form abbreviations.

The key outcomes of the present study are that in the considered test cases:
\begin{itemize}
\item W-SIA allows larger time steps compared to standard SIA but retains quadratic scaling,
\item W-SIA is slightly more accurate than SIA,
\item W-SIA-FSSA allows for significantly larger time steps with linear scaling limited to coarse $\Delta x$.
\end{itemize}
Furthermore we find that:
\begin{itemize}
\item W-SIAStokes has a linear time step restriction scaling, but now for all $\Delta x$ (no FSSA),
\item W-SIAStokes-FSSA has a linear time step restriction scaling for all $\Delta x$,
\item W-SIAStokes is numerically more robust compared to W-SIA for a negligible cost increase,
\item W-SIAStokes is more accurate than W-SIA.
\end{itemize}
We expand on the two lists above in the paragraphs written below. 

The first weak formulation is W-SIA \eqref{eq:sia_weakform_less}. 
An immediate benefit when using W-SIA is that we were able to add the FSSA stabilization terms \eqref{eq:fssa} (W-SIA-FSSA). 
This improved the time step restriction from quadratic (W-SIA) to linear (W-SIA-FSSA) scaling in terms of the horizontal mesh size $\Delta x$ 
when using the W-SIA velocities for solving the free surface equation. 
Overall, the time step sizes in all the considered test cases were increased for at least approximately $100$-times 
as compared to the standard SIA formulation time step sizes.. 
However, we observed that when $\Delta x$ is small enough, the largest time step size behavior became unpredictable 
and started taking values similar to the standard SIA time step sizes. 

As a remedy, we extended W-SIA to W-SIAStokes \eqref{eq:sia_weakform_all} 
by adding the originally neglected stress terms back to W-SIA, but kept 
the SIA viscosity \eqref{eq:SIAvisco} intact so that W-SIAStokes remained a linear problem. 
We observed predictable largest time step size behavior, which had a linear scaling in terms of $\Delta x$ without even adding the FSSA 
stabilization terms to the weak formulation. 
After we added the FSSA stabilization terms to W-SIAStokes (W-SIAStokes-FSSA) the scaling remained linear, and 
the largest time step sizes increased for at least approximately $100$-times as compared to the standard SIA 
formulation time step sizes. In one of the tests we compared W-SIAStokes to W-Stokes and found that the time step sizes, 
including the scaling in terms of $\Delta x$, were comparable.

When compared to the standard SIA model, W-SIA and W-SIAStokes 
models had a smaller error, taking the W-Stokes solution as reference. This held true also when the FSSA stabilization terms were added 
to the two weak formulations and very long time steps were used. The error vs. runtime ratio was also favourable in the 
case of both weak formulations (with and without the FSSA stabilization) over the standard SIA solution. 
Among the two, W-SIAStokes had a better error vs. runtime ratio, 
(with and without the FSSA stabilization). We note that the runtime measurements could differ depending on the 
code implementation of the different SIA formulations. 

To view the observed runtime measurements from another angle, we performed a theoretical computational cost estimation 
in Section \ref{section:performance_analysis}. Based on that we conclude that the computational cost 
of W-SIA and W-SIAStokes is comparable to that of the standard SIA model, 
in terms of the asymptotic behaviour, that is, when the number of the horizontal mesh vertices is increased.

We anticipate that glaciologists interested in the ice spin-up simulations of which the simulation length is on 
the order of hundreds of thousands of years, 
could benefit from using the W-SIA-FSSA model or the W-SIAStokes-FSSA model over the standard SIA model. 
Note that our results are for isothermal simulations. More studies are needed to investigate the time step restrictions 
related to the temperature evolution and other physical processes. However we believe resolving the restriction due 
to the velocity-surface coupling is the most difficult problem.

We speculate that the W-SIAStokes model is a good choice to replace the standard SIA model within the scope of coupled models. 
One of them is the ISCAL (Ice Sheet Coupled Approximation Level) model \cite{iscal}
where the standard SIA model and the W-Stokes model are used dynamically over an ice sheet geometry, depending on 
the desired modeling accuracy. When both models are coupled together, the time step restriction is bound to that of the SIA model 
(quadratic scaling in terms of $\Delta x$), whereas the W-SIAStokes model allows for linear scaling.



\section*{Acknowledgments}
Josefin Ahlkrona and Igor Tominec were funded by the the Swedish Research Council, grant number 2021-04001. 
Josefin Ahlkrona also received funding from Swedish e-Science Research Centre (SeRC).

\appendix
\section{Von Neumann Stability Analysis}
\label{sec:appendix:von_neumann}
\subsection{Without FSSA}
\subsubsection{Coupling}

The solution to the strong form shallow-ice approximation (SIA) 

\begin{alignat}{2}
- \frac{\partial p}{\partial x}+\frac{\partial }{\partial y}\left(\mu\frac{\partial u_x}{\partial y}\right)&=0 \label{eq:SIAeqa} \\
-\frac{\partial p}{\partial y} &= \rho g \label{eq:SIAeqb}
\end{alignat}
is in the isothermal case
\begin{align*}
p &=\rho g (h-y)\\
u_1 &= - \frac{1}{2}\mathfrak{A} (\rho g)^3 \left(\frac{\partial h}{ \partial x}\right)^3  \left( (h-b)^4- (h-y)^4 \right),\\
u_2&=\int_b^y \frac{\partial u_x}{\partial x} dy = \frac{1}{2}\mathfrak{A} (\rho g)^3\left(3 \left( \frac{\partial h}{\partial x} \right)^2 \frac{\partial^2 h}{\partial x^2} \left( (h-b)^4 (y-b)+\frac{(h-y)^5}{5}-\frac{(h-b)^5}{5}\right) \right. \\
& + \left.4 \left( \frac{\partial h}{\partial x} \right)^3 \left((h-b)^3\left( \frac{\partial h}{\partial x}-\frac{\partial b}{\partial x} \right) (y-b) + \frac{(h-y)^4}{4}\frac{\partial h}{\partial x}-\frac{(h-b)^4}{4}\frac{\partial h}{\partial x}\right)
\right)
\end{align*}
defining $H=(h-b)$ and cosindering only the surface $y=h$ we get:
\begin{align*}
u_1 &=- \frac{1}{2}\mathfrak{A} (\rho g)^3 \left(\frac{\partial h}{ \partial x}\right)^3   H^4 \\
u_2&=\frac{1}{2}\mathfrak{A} (\rho g)^3\left(3 \left( \frac{\partial h}{\partial x} \right)^2 \frac{\partial^2 h}{\partial x^2} \left( H^5 -\frac{H^5}{5}\right) +4 \left( \frac{\partial h}{\partial x} \right)^3 \left(H^4 \frac{\partial H}{\partial x}  -\frac{H^4}{4}\frac{\partial h}{\partial x}\right)\right)
\end{align*}
Inserting the closed form expressions for the into the time discretization of the free-surface equation \eqref{eq:fssa_free_surface_time_discretized} revals the full coupled system
\begin{alignat}{2}
&\frac{h^{k+1}-h^k}{\Delta t}+\left(-\frac{1}{2}A(\rho g)^3 \left(\frac{\partial h}{\partial x}\right)^3  (h-b)^4\right) \frac{\partial h^{k+\gamma}}{\partial x} \\
&=\frac{1}{2}\mathfrak{A} (\rho g)^3\left(3 \left( \frac{\partial h}{\partial x} \right)^2 \frac{\partial^2 h}{\partial x^2} \left( H^5 -\frac{H^5}{5}\right) +4 \left( \frac{\partial h}{\partial x} \right)^3 \left(H^4 \frac{\partial H}{\partial x}  -\frac{H^4}{4}\frac{\partial h}{\partial x}\right)\right) + a_s
\label{eq:fullsystem}
\end{alignat}
This is a highly non-linear equation and needs to be linearized before Fourier analysis can be used. 

Wihtout FSSA all $h$ and $H$-terms which do not have a superscript $(\cdot)^{k+\gamma}$ are approximated explicitely, i.e. they all will have a superscipt $(\cdot)^{k}$.

\subsubsection{Linearization}

Following \citep{CHENG201729} (where the ice thickness equation was analyzed) we consider a slab on a slope with a small 
surface perturbation $\delta$ around an average state $\bar{h}$ which is just an inclined plane. We write:

\begin{alignat}{2}
h&=\bar{h}+\delta \text{ (for $h$ related to velocity)}, \quad  H=\bar{H}+\delta \\
h&=\bar{h}+\hat{\delta} \text{ (for $h$ explicitely in the free surface equation)} , \quad  H=\bar{h}+\hat{\delta }.
\end{alignat}
Inserting in \eqref{eq:fullsystem} yields
\begin{align*}
&\frac{\partial (\bar{h}+\hat{\delta})}{\partial t}+\left(-\frac{1}{2}A(\rho g)^3 \left(\frac{\partial (\bar{h}+\delta)}{\partial x}\right)^3  ((\bar{H}+\delta))^4\right) \frac{\partial (\bar{h}+\hat{\delta})}{\partial x} \\
&= \frac{1}{2}\mathfrak{A} (\rho g)^3\left(3 \left( \frac{\partial \bar{h}+\delta}{\partial x} \right)^2 \frac{\partial^2 \bar{h}+\delta}{\partial x^2} \left( (\bar{H}+\delta)^5 -\frac{(\bar{H}+\delta)^5}{5}\right) \right.\\ 
&\left.+4 \left( \frac{\partial \bar{h}+\delta}{\partial x} \right)^3 \left((\bar{H}+\delta)^4 \frac{\partial (\bar{H}+\delta)}{\partial x}  -\frac{(\bar{H}+\delta)^4}{4}\frac{\partial \bar{h}+\delta  }{\partial x}\right)\right) + a_s
\end{align*}
Using that the second derivative of the steady state surface $\bar{h}$ is zero we get
\begin{align*}
&\frac{\partial (\hat{\delta})}{\partial t}+\left(-\frac{1}{2}A(\rho g)^3 \left(\frac{\partial \bar{h}}{\partial x}\right)^3  (\bar{H})^4\right) \frac{\partial \hat{\delta}}{\partial x} 
+\left(-\frac{1}{2}A(\rho g)^3 \left(\frac{\partial (\bar{h}+\delta)}{\partial x}\right)^3  (\bar{H})^4\right) \frac{\partial \bar{h}}{\partial x} \\
&+\left(-\frac{1}{2}A(\rho g)^3 \left(\frac{\partial \bar{h}}{\partial x}\right)^3  ((\bar{H}+\delta))^4\right) \frac{\partial \bar{h}}{\partial x} \\
&= \frac{1}{2}\mathfrak{A} (\rho g)^3\left(3 \left( \frac{\partial \bar{h}}{\partial x} \right)^2 \frac{\partial^2 \delta}{\partial x^2}  \frac{4}{5}\bar{H}^5 +4 \left( \frac{\partial \bar{h}+\delta}{\partial x} \right)^3 \left((\bar{H})^4 \frac{\partial (\bar{H})}{\partial x}  -\frac{(\bar{H})^4}{4}\frac{\partial \bar{h} }{\partial x}\right)\right) \\
&\left.+4 \left( \frac{\partial \bar{h}}{\partial x} \right)^3 \left((\bar{H}+\delta)^4 \frac{\partial (\bar{H})}{\partial x}+(\bar{H})^4 \frac{\partial (\bar{H}+\delta)}{\partial x}  -\frac{(\bar{H}+\delta)^4}{4}\frac{\partial \bar{h} }{\partial x}-\frac{(\bar{H})^4}{4}\frac{\partial \bar{h}+\delta }{\partial x}\right)\right) + a_s
\end{align*}

And ignoring higher order terms in $\delta$ furthermore yields

\begin{align*}
&\frac{\partial (\hat{\delta})}{\partial t}+\left(-\frac{1}{2}A(\rho g)^3 \left(\frac{\partial \bar{h}}{\partial x}\right)^3  (\bar{H})^4\right) \frac{\partial \hat{\delta}}{\partial x} 
+\left(-\frac{1}{2}A(\rho g)^3 3\left(\frac{\partial (\bar{h})}{\partial x}\right)^2\left(\frac{\partial (\delta)}{\partial x}\right)  (\bar{H})^4\right) \frac{\partial \bar{h}}{\partial x} 
+\left(-\frac{1}{2}A(\rho g)^3 \left(\frac{\partial \bar{h}}{\partial x}\right)^3  4\bar{H}^3\delta \right) \frac{\partial \bar{h}}{\partial x} \\
&= \frac{1}{2}\mathfrak{A} (\rho g)^3\left(3 \left( \frac{\partial \bar{h}}{\partial x} \right)^2 \frac{\partial^2 \delta}{\partial x^2}  \frac{4}{5}\bar{H}^5 +12 \left( \frac{\partial \bar{h}}{\partial x} \right)^2 \frac{\partial \delta}{\partial x} \left((\bar{H})^4 \frac{\partial (\bar{H})}{\partial x}  -\frac{(\bar{H})^4}{4}\frac{\partial \bar{h} }{\partial x}\right)\right. \\
&\left.+4 \left( \frac{\partial \bar{h}}{\partial x} \right)^3 \left(4\bar{H}^3\delta \frac{\partial (\bar{H})}{\partial x}+(\bar{H})^4 \frac{\partial (\delta)}{\partial x}  -\bar{H}^3\delta \frac{\partial \bar{h} }{\partial x}-\frac{(\bar{H})^4}{4}\frac{\partial \delta }{\partial x}\right)\right) + a_s
\end{align*}
Considering that the steady state thickness is zero further simplifies the expression into
\begin{align*}
&\frac{\partial \hat{\delta}}{\partial t}+\left(-\frac{1}{2}A(\rho g)^3 \left(\frac{\partial \bar{h}}{\partial x}\right)^3  \bar{H}^4\right) \frac{\partial \hat{\delta}}{\partial x} +\left(-\frac{1}{2}A(\rho g)^3 3\left(\frac{\partial \bar{h}}{\partial x}\right)^2\left(\frac{\partial (\delta)}{\partial x}\right)  \bar{H}^4\right) \frac{\partial \bar{h}}{\partial x} +\left(-\frac{1}{2}A(\rho g)^3 \left(\frac{\partial \bar{h}}{\partial x}\right)^3  4\bar{H}^3\delta \right) \frac{\partial \bar{h}}{\partial x} \\
&= \frac{1}{2}\mathfrak{A} (\rho g)^3\left(3 \left( \frac{\partial \bar{h}}{\partial x} \right)^2 \frac{\partial^2 \delta}{\partial x^2}  \frac{4}{5}\bar{H}^5 -3 \left( \frac{\partial \bar{h}}{\partial x} \right)^2 \frac{\partial \delta}{\partial x}   \bar{H}^4\frac{\partial \bar{h} }{\partial x}+4 \left( \frac{\partial \bar{h}}{\partial x} \right)^3 \left(\frac{3}{4}\bar{H}^4 \frac{\partial \delta}{\partial x}  -\bar{H}^3\delta \frac{\partial \bar{h} }{\partial x}\right)\right) + a_s
\end{align*}

rearranging and defining $\frac{\partial \bar{h}}{\partial x}=C_\alpha$ we get:

\begin{align*}
&\frac{\partial \hat{\delta}}{\partial t}-\frac{1}{2}\mathfrak{A}(\rho g)^3 C_\alpha^3  \bar{H}^4 \frac{\partial \hat{\delta}}{\partial x} -\frac{3}{2}\mathfrak{A}(\rho g)^3 C_\alpha^3  \bar{H}^4 \frac{\partial \delta}{\partial x} = \frac{6}{5}\mathfrak{A} (\rho g)^3  C_\alpha^2  \bar{H}^5 \frac{\partial^2 \delta}{\partial x^2}   + a_s
\end{align*}
which we will write on the form
\begin{align*}\label{eq:linearizedeq}
&\frac{\partial \hat{\delta}}{\partial t}-C_1 \frac{\partial \hat{\delta}}{\partial x} -C_2 \frac{\partial \delta}{\partial x} = C_3\frac{\partial^2 \delta}{\partial x^2}   + a_s
\end{align*}

with $C_1=\frac{1}{2}\mathfrak{A}(\rho g)^3 C_\alpha^3  \bar{H}^4$, $C_2= \frac{3}{2}\mathfrak{A}(\rho g)^3 C_\alpha^3  \bar{H}^4$ and $C_3=\frac{6}{5}\mathfrak{A} (\rho g)^3  C_\alpha^2  \bar{H}^5 $. OBS! $C_1$ and $C_2$ are negative!

We can discretize $\hat{\delta}$ using $k$ or $k+1$ (the latter is better of course), while $\delta$  must be taken  from time-step $k$ as it originates from the velocity. 

\subsubsection{Fourier analysis}
We will now consider the SIA solution in Fourier space. We will thus apply a  Fourier transform, and consider one frequency at a time 
\begin{align*}
\delta^k_j \rightarrow \tilde{\delta}^k e^{i n x_j}=\tilde{\delta} e^{i n  j \Delta x}\\
\delta^k_{j+1} \rightarrow \tilde{\delta}^k e^{i n x_{j+1}}=\tilde{\delta}^k e^{i n  j \Delta x}e^{i n  \Delta x}\\
\end{align*}
The factor $e^{i n  j \Delta x}$ will appear in every term of every equation, and we will hencedivide by that and not write it out from here on. We will also consider the finite difference discretization which best corresponds to a P1 FEM discretization, meaning that second derivatives will be represented as central differences, so that $\frac{\partial \delta}{\partial x}$ and $\frac{\partial^2 \delta}{\partial x^2} $ are discretized and transformed as:
\begin{align*}
\frac{\partial \delta}{\partial x} \approx \frac{\delta_{j+1}-\delta_{j-1}}{2 \Delta x}  \rightarrow \tilde{\delta}\frac{e^{in \Delta x}-e^{-in \Delta x}}{2 \Delta x} = \tilde{\delta}\frac{2i sin (n\Delta x)}{2 \Delta x}\\
\frac{\partial^2 \delta}{\partial x^2}\approx \frac{\delta_{j+1}-2\delta_{j}+\delta_{j-1}}{(\Delta x)^2}  \rightarrow \tilde{\delta} \frac{  e^{in \Delta x}-2+e^{-in \Delta x} }{(\Delta x)^2} = -\tilde{\delta}\frac{4 \sin^2 (n \Delta x/2)}{ (\Delta x)^2  }
\end{align*}
Inserting this into the linearized equation equation \label{eq:linearizedeq}, assuming implicit handling of the free surface equation $\alpha=1$ itself and explicit handling of velocities, and setting $a_s=0$ gives

\begin{align*}
&\frac{\tilde{\delta}^{k+1}-\tilde{\delta}^{k}}{\Delta t}-C_1  \tilde{\delta}^{k+1}\frac{2i sin (n\Delta x)}{2 \Delta x}      -C_2  \tilde{\delta}^{k}\frac{2i sin (n\Delta x)}{2 \Delta x}  =- C_3 \tilde{\delta^k}\frac{4 \sin^2 (n \Delta x/2)}{ (\Delta x)^2  }.
\end{align*}

Using Euler's formulas and rearranging gives:


\begin{align*}
&\tilde{\delta}^{k+1} \left(1 -\Delta t C_1 \frac{2i sin (n\Delta x)}{2 \Delta x} \right)      =\tilde{\delta}^{k}  \left(1 - \Delta tC_3 \frac{4 \sin^2 (n \Delta x/2)}{ (\Delta x)^2  }+ \Delta t C_2 \frac{2i sin (n\Delta x)}{2 \Delta x}  \right)   .
\end{align*}
In order for each Fourier mode to stay bounded as the simulation runs, i.e. that $|\tilde{\delta}^{k+1}|\leq |\tilde{\delta}^{k}|$, we require that
\begin{equation}
\frac{\left|1 - \Delta tC_3 \frac{4 \sin^2 (n \Delta x/2)}{ (\Delta x)^2  }+ \Delta t C_2 \frac{2i sin (n\Delta x)}{2 \Delta x}  \right| }{ \left|1 -\Delta t C_1 \frac{2i sin (n\Delta x)}{2 \Delta x} \right|} \leq 1.
\end{equation}
The denominator value of the parenthesis of the left hand side is always larger than one (which is due to that $\alpha=1$). We thus oly need to boundthe nominator. The complex part of the nominator is bounded  as long as    

\begin{align*}
\frac{|C_2| \Delta t }{2\Delta t}\leq  1,
\end{align*}
while the real part needs to fulfill 
\begin{align*}
\left|1 - \Delta tC_3 \frac{4 \sin^2 (n \Delta x/2)}{ (\Delta x)^2  }\right| <1 \Rightarrow \frac{C_3\Delta t}{(\Delta x)^2}<0.5.
\end{align*}
So depending on the balance between $C_2$ and $C_3$ we get a linear or parabolic time-step constraint. $C_3$ is big for thick ice with flat slopes, i.e. for the dynamics typical of the interior of an ice sheet, we get a parabolic constraint.

\subsection{With FSSA}
\subsubsection{Coupling}

As observed from \eqref{eq:FSSAmakesPimplicit_step2} the FSSA make the discretization of the pressure implicit, i.e.
For $\theta =1$ we have 
\begin{equation}
p(x,y,t) \approx \rho g (h^{k+1} - y).
\end{equation}
The pressure yields $u_1$ and $u_1$ gives $u_2$, so this implicit discretiation will propagate to the velocity solution in the following way
\begin{align*}
u_1 &= - \frac{1}{2}\mathfrak{A} (\rho g)^3 \left(\frac{\partial h^{k+1}}{ \partial x}\right)^3  \left( (h^k-b)^4- (h^k-y)^4 \right),\\
%
%
%
u_2 & =\int_b^y \frac{\partial u_x}{\partial x} dy = 3 \left( \frac{\partial h^{k+1}}{\partial x} \right)^2 \frac{\partial^2 h^{k+1}}{\partial x^2} \left( (h^k-b)^4 (y-b)+\frac{(h^k-y)^5}{5}-\frac{(h^k-b)^5}{5}\right)  \\
&  +4 \left( \frac{\partial h^{k+1}}{\partial x} \right)^3 \left((h^k-b)^3\left( \frac{\partial h^k}{\partial x}-\frac{\partial b}{\partial x} \right) (y-b) + \frac{(h^k-y)^4}{4}\frac{\partial h^k}{\partial x}-\frac{(h^k-b)^4}{4}\frac{\partial h^k}{\partial x}\right)
%
\end{align*}

Then the final expressions at $y=h$ are :

\begin{align*}
u_1 &=- \frac{1}{2}\mathfrak{A} (\rho g)^3 \left(\frac{\partial h^{k+1}}{ \partial x}\right)^3   (H^k)^4 \\
u_2&=\frac{1}{2}\mathfrak{A} (\rho g)^3\left(3 \left( \frac{\partial h^{k+1}}{\partial x} \right)^2 \frac{\partial^2 h^{k+1}}{\partial x^2} \left( (H^k)^5 -\frac{(H^k)^5}{5}\right) +4 \left( \frac{\partial h^{k+1}}{\partial x} \right)^3 \left((H^k)^4 \frac{\partial H^k}{\partial x}  -\frac{(H^k)^4}{4}\frac{\partial h^k}{\partial x}\right)\right)
\end{align*}
and the full coupled system is then
\begin{alignat}{2}
&\frac{h^{k+1}-h^k}{\Delta t}+\left(-\frac{1}{2}A(\rho g)^3 \left(\frac{\partial h^{k+1}}{\partial x}\right)^3  (H^k)^4\right) \frac{\partial h^{k+\gamma}}{\partial x} \\
&=\frac{1}{2}\mathfrak{A} (\rho g)^3\left(3 \left( \frac{\partial h^{k+1}}{\partial x} \right)^2 \frac{\partial^2 h^{k+1}}{\partial x^2} \left( (H^k)^5 -\frac{H^5}{5}\right) +4 \left( \frac{\partial h^{k+1}}{\partial x} \right)^3 \left((H^k)^4 \frac{\partial H^k}{\partial x}  -\frac{(H^k)^4}{4}\frac{\partial h^k}{\partial x}\right)\right) + a_s
\label{eq:fullsystemwFSSA}
\end{alignat}

\subsubsection{Linearization}

With the same approach as before we get
\begin{align*}
&\frac{\partial \hat{\delta}}{\partial t}+\left(-\frac{1}{2}A(\rho g)^3 \left(\frac{\partial \bar{h}}{\partial x}\right)^3  \bar{H}^4\right) \frac{\partial \hat{\delta}}{\partial x} \\
&+\left(-\frac{1}{2}A(\rho g)^3 3\left(\frac{\partial \bar{h}}{\partial x}\right)^2\left(\frac{\partial \delta^{k+1}}{\partial x}\right)  \bar{H}^4\right) \frac{\partial \bar{h}}{\partial x} \\
&+\left(-\frac{1}{2}A(\rho g)^3 \left(\frac{\partial \bar{h}}{\partial x}\right)^3  4\bar{H}^3\delta^k \right) \frac{\partial \bar{h}}{\partial x} \\
&= \frac{1}{2}\mathfrak{A} (\rho g)^3\left(3 \left( \frac{\partial \bar{h}}{\partial x} \right)^2 \frac{\partial^2 \delta^{k+1}}{\partial x^2}  \frac{4}{5}\bar{H}^5 -3 \left( \frac{\partial \bar{h}}{\partial x} \right)^2 \frac{\partial \delta^{k+1}}{\partial x}   \bar{H}^4\frac{\partial \bar{h} }{\partial x}+4 \left( \frac{\partial \bar{h}}{\partial x} \right)^3 \left(\frac{3}{4}\bar{H}^4 \frac{\partial \delta^k}{\partial x}  -\bar{H}^3\delta^k \frac{\partial \bar{h} }{\partial x}\right)\right) + a_s
\end{align*}
rearranging and setting $\frac{\partial \bar{h}}{\partial x}=\frac{\partial \bar{h}}{\partial x}=C_\alpha$ yields
\begin{align*}
&\frac{\partial \hat{\delta}}{\partial t}-\frac{1}{2}\mathfrak{A}(\rho g)^3 C_\alpha^3  \bar{H}^4 \frac{\partial \hat{\delta}}{\partial x} -\frac{3}{2}\mathfrak{A}(\rho g)^3 C_\alpha^3  \bar{H}^4 \frac{\partial \delta^{k}}{\partial x} = \frac{6}{5}\mathfrak{A} (\rho g)^3  C_\alpha^2  \bar{H}^5 \frac{\partial^2 \delta^{k+1}}{\partial x^2}   + a_s
\end{align*}
which we will write on the form
\begin{align*}
&\frac{\partial \hat{\delta}}{\partial t}-C_1 \frac{\partial \hat{\delta^{k+1}}}{\partial x} -C_2 \frac{\partial \delta}{\partial x} = C_3\frac{\partial^2 \delta^{k+1}}{\partial x^2}   + a_s
\end{align*}
So the difference in the linearized equation is that the second derivative is treated implicitly. 

\subsubsection{Fourier analysis}
We get the same expression as without FSSA, only that the second derivative is now treated implicitly
\begin{align*}
&\frac{\tilde{\delta}^{k+1}-\tilde{\delta}^{k}}{\Delta t}-C_1  \tilde{\delta}^{k+1}\frac{2i sin (n\Delta x)}{2 \Delta x}      -C_2  \tilde{\delta}^{k}\frac{2i sin (n\Delta x)}{2 \Delta x}  =- C_3 \tilde{\delta^k}\frac{4 \sin^2 (n \Delta x/2)}{ (\Delta x)^2  }.
\end{align*}
Rearranging gives
%
\begin{align*}
&\tilde{\delta}^{k+1} \left(1 -\Delta t C_1 \frac{2i sin (n\Delta x)}{2 \Delta x}+ \Delta tC_3 \frac{4 \sin^2 (n \Delta x/2)}{ (\Delta x)^2  } \right)      =\tilde{\delta}^{k}  \left(1 + \Delta t C_2 \frac{2i sin (n\Delta x)}{2 \Delta x}  \right)   
\end{align*}
The absolute value of the parenthesis to the left is clearly negative, and the one on the right hand side has an absolute value smaller than one if $\frac{\Delta t}{\Delta x} C_2\leq 1\Rightarrow \Delta t \leq \frac{\Delta x}{ \frac{3}{2}\mathfrak{A}(\rho g)^3 C_\alpha^3  \bar{H}^4}$, i.e. we get a linear time-step constraint.

\bibliographystyle{abbrvnat}
\bibliography{sections/references_igorspaper}

\end{document}